\documentclass{report}
\usepackage{amsmath}
\usepackage{amsfonts}
\usepackage{amsthm}
\usepackage{amssymb}

\usepackage{epsfig}

\usepackage{makeidx}
\makeindex

\theoremstyle{definition}
\newtheorem{definition}{Definition}[chapter]

\newtheorem*{example}{Example}

\theoremstyle{plain}
\newtheorem{lemma}[definition]{Lemma}
\newtheorem{theorem}[definition]{Theorem}
\newtheorem{corollary}[definition]{Corollary}

\theoremstyle{remark}
\newtheorem*{remark}{Remark}

\newcounter{listcounter}
\newenvironment{myenumerate}
 {\begin{list}{\textnormal{(\thelistcounter)}}{\usecounter{listcounter}}}
 {\end{list}}

\newcommand{\fix}{\mathrm{fix}\,}
\newcommand{\fixw}{\mathrm{fix}_w\,}
\newcommand{\ch}{\mathrm{ch}}
\newcommand{\pz}[1]{\frac{p_{#1}}{z_{#1}}}
\newcommand{\pzm}[2]{\frac{p_{#2} (#1)}{z_{#2}}}
\newcommand{\trace}{\mathrm{tr}\,}
\newcommand{\spec}[1]{\mathbf{S}_{#1}}
\newcommand{\mspec}[1]{\mathfrak{M}_{#1}}
\newcommand{\vspec}[1]{\mathbf{V}_{#1}}
\newcommand{\specf}[1]{\mathbf{S}^{*}_{#1}}
\newcommand{\vspecf}[1]{\mathbf{V}^{*}_{#1}}
\newcommand{\partitions}{\mathcal{P}}
\newcommand{\genlin}{\operatorname{GL}}
\newcommand{\repring}{\operatorname{R}}
\newcommand{\centfunc}[2]{\operatorname{CF}_{#2}(#1)}
\newcommand{\charring}[1]{\operatorname{CR}(#1)}
\newcommand{\inpleth}[1]{\,{\boxtimes}_{#1}}
\newcommand{\inv}[1]{{#1}^{-1}}
\newcommand{\orbit}[2]{\mathcal{O}_{#1}(#2)}
\newcommand{\stab}[2]{\mathrm{stab}_{#1}(#2)}
\newcommand{\catsets}[1]{\mathrm{Sets}_{#1}}
\newcommand{\catfsets}[1]{\mathbb{B}_{#1}}
\newcommand{\posints}{\mathbb{P}}
\newcommand{\cmplx}{\mathbb{C}}
\newcommand{\restrict}[2]{{#1}{|}_{#2}}
\newcommand{\myprime}[1]{{#1}'}
\newcommand{\cardw}[1]{|{#1}{|}_{w}}
\newcommand{\adams}[1]{{\psi}^{#1}}
\newcommand{\pp}[3]{{p_{#1}(#2)}^{#3}}
\newcommand{\digraphs}[1]{D_{#1}}
\newcommand{\digraphsl}[1]{D^{\circ}_{#1}}
\newcommand{\aut}[1]{\operatorname{Aut}(#1)}
\newcommand{\graphs}[1]{\mathcal{G}_{#1}}
\newcommand{\graphsl}[1]{\mathcal{G}^{\circ}_{#1}}
\newcommand{\bigraphs}[1]{\mathcal{B}_{#1}}

\newcommand{\lm}{\ensuremath{\lambda}}

\newcommand{\sym}[1]{\ensuremath{{\mathfrak{S}_{#1}}}}
\newcommand{\lam}[2]{\ensuremath{\Lambda^{#1}_{#2}}}
\newcommand{\lamhat}[1]{\ensuremath{\hat{\Lambda}_{#1}}}

\newcommand{\lamring}{$\lambda$-ring}
\newcommand{\lamrings}{$\lambda$-rings}
\newcommand{\myiff}{if and only if}
\newcommand{\maple}{Maple}
\newcommand{\mcc}[1]{\multicolumn{1}{c|}{#1}}

\begin{document}

	\newpage
 \renewcommand{\thepage}{}

	\pagestyle{empty}
	\begin{center}
	\vspace*{.6in}
	   Graphical Enumeration: A Species-Theoretic Approach \\
	\vspace*{.75in}
	   A Dissertation \\
	\vspace*{.5in}
           Presented to \\
	\vspace*{.1in}
           The Faculty of the Graduate School of Arts and Sciences \\
	\vspace*{.1in}
           Brandeis University \\
	\vspace*{.1in}
           Department of Mathematics \\
	\vspace*{.1in}
          Ira Gessel, Advisor \\
	\vspace*{.75in}
	   In Partial Fulfillment \\
	\vspace*{.05in}
	   of the Requirements for the Degree\\
	\vspace*{.05in}
	   Doctor of Philosophy\\
	\vspace*{.75in}
	   by\\
	\vspace*{.1in}
	   Leopold Travis\\
	\vspace*{.1in}
	   February, 1999\\
	\end{center} 

\newpage
\pagestyle{plain}
\setcounter{page}{3}
\renewcommand{\thepage}{\roman{page}}

	\begin{center}
	   \vspace*{.25in}
             {\large \bf ACKNOWLEDGEMENTS }
	\end{center}
     
\vspace*{.25in}

Above all, I would like to thank my advisor, Ira Gessel, for his
patience and encouragement through some difficult times
--- and for sharing so many good ideas with me.  
Without him, this work would never have taken shape.

\newpage
	\begin{center}
	   \vspace*{.25in}
             {\large \bf ABSTRACT }\\
	   \vspace*{.25in}
             {\bf  Graphical Enumeration: A Species-Theoretic Approach } \\
	   \vspace*{.25in}
	\end{center}

	\begin{center}
          \noindent
          \begin{minipage}{.85\linewidth}
          \raggedright
            A dissertation presented to the Faculty of
            the Graduate School of Arts and Sciences of
            Brandeis University, Waltham, Massachusetts
          \end{minipage}
	\end{center}

\vspace*{.25in}

	\begin{center}
            by Leopold Travis
	\end{center}

\vspace*{.25in}

An operation on species corresponding to the inner plethysm
of their associated cycle index series is constructed.
This operation, the inner plethysm of species, is generalized to
$n$-sorted species.
Polynomial maps on species are studied and used to extend
inner plethysm and other operations to virtual species.
Finally, inner plethysm and other operations on species are
applied to various problems in graph theory.  

In particular,
regular graphs, and digraphs in which every vertex has outdegree $k$,
are enumerated.

\tableofcontents

\listoftables

\listoffigures

\chapter{Introduction}

\setcounter{page}{1}
\renewcommand{\thepage}{\arabic{page}}


Since its introduction by Joyal in \cite{spec1}, the concept of
\emph{species of structures} has proved useful in many areas of
combinatorics. Our aim in this work is two-fold: to study certain new or
little-known operations in the theory of species, and to apply these operations
to various problems in graphical enumeration.

Chapter~\ref{chap:spec} is devoted to species. We
define the cycle index series (in a slightly different way from Joyal), 
and relate it to the \emph{analytic functors} he introduces
in~\cite{anfunc}. We use our techniques to prove the well-known formula
for the cycle index series of the composition of two species.

We then investigate the operation of \emph{inner plethysm}. As an 
operation on symmetric functions it was introduced by 
Littlewood in \cite{littlewood}, and has been used in the study of representation
theory (see Thibon~\cite{thibon}). The corresponding operation on
species has not been examined. We give two constructions for it,
one using combinatorial operations on species, the other using analytic
functors.

Next, we examine inner plethysm in the context
of $n$-sorted species.  We focus upon the
analytic functor approach, using it both to define the operation of
\emph{inner plethysm in $Y$}, and to calculate the corresponding operation
on symmetric functions.

Finally, we examine \emph{polynomial maps on species}, using
essentially the techniques of \cite{virtspec} to extend our results
to virtual species.

In Chapter~\ref{chap:digraphs} we 
study of digraphs. In particular we count unlabeled digraphs in which
every vertex has outdegree $k$, which to our knowledge is an open problem.

In Chapter~\ref{chap:graphs} we study graphs, once again emphasizing
a species-theoretic point of view and applying techniques from 
Chapter~\ref{chap:spec}.

\chapter{Species}
\label{chap:spec}

\section{The cycle index and symmetric functions}
\label{sec:cycsym}

Following \cite{symfunc}, we denote by $\lam{}{}$ the ring of symmetric functions
in the infinite set of indeterminates $x_1,x_2,\dots$, and
by $\lam{n}{}$ the subring of $\lam{}{}$ consisting of symmetric 
functions homogeneous of degree $n$. More generally, $\lam{}{x}$
will denote the ring of symmetric functions in the indeterminates
$x_1,x_2,\dots$, $\lam{}{y}$  will denote the ring of symmetric
functions in $y_1,y_2,\dots$, etc. We will denote by $\lam{}{xy}$ the ring of 
formal series in $x_1,x_2,\dots,y_1,y_2,\dots$ of bounded degree, which are
symmetric separately in the $x_i$ and in the $y_i$. We define
$\lam{}{xyz}$, etc., similarly.
It is not difficult to show that
$\lam{}{xy}$ is (isomorphic to) the tensor product
$\lam{}{x} \otimes \lam{}{y}$, $\lam{}{xyz} \cong
\lam{}{x} \otimes \lam{}{y} \otimes \lam{}{z}$, and so forth.
We define $\lamhat{}$ to be the ring of symmetric functions
in $x_1,x_2,\dots$ of unbounded degree, and $\lamhat{xy}$, etc.,
similarly.

Let $F$ be a species, that is, a functor from $\catfsets{}$ to $\catfsets{}$,
where $\catfsets{}$ is the category of finite sets and bijections.
We denote by $[n]$ the set $\{1,\dots,n\}$, and
write $F[n]$ for $F[\{1,\dots,n\}]$. The symmetric
group $\sym{n}$ acts on $F[n]$, since by functoriality any $\sigma \in
\sym{n}$ induces a permutation $F[\sigma]$ of $F[n]$. For any partition
$\lm$ of $n$, we define $\fix F[\lm]$
to be the number of fixed points of $F[\sigma]$, where $\sigma$ is
any permutation of cycle type $\lm$. Following~\cite{lagrangeinv}, we associate
to $F$ its \emph{cycle index series}, or
\emph{cycle index}, the symmetric function $Z_F$:

\begin{equation}
Z_F = \sum_{\lm} \fix F[\lm] \pz{\lm}
\label{eq:cycdef}
\end{equation}
In this and all similar sums, the index of summation $\lm$ is taken
to run over the set $\partitions$ of all partitions.

Our definition of the cycle index differs from that introduced by Joyal
in~\cite{spec1} in that we use the power sum symmetric function $p_i$
in place of an indeterminate. 
Thus our cycle index series can be
considered both as a formal series in the power sums $p_i$, and
as a symmetric function in the underlying variables $x_i$. 
We will explore the combinatorial significance of this second interpretation
in the next sections. First, we give a slight generalization of~\eqref{eq:cycdef},
to
weighted species (see~\cite[Section 6]{spec1}).

For a
commutative ring $R$ containing the rationals, let $\catsets{R}$ denote
the category of \emph{$R$-weighted sets}. The objects of $\catsets{R}$
are pairs $(A,w)$, where $A$ is a set and $w:A\to R$ is a function, which
we will refer to as a \emph{weight function}. We do not require that $A$
be finite, but $\sum_{a\in A} w(a)$ must exist in $R$ for
$(A,w)$ to be an object of $\catsets{R}$; we denote this sum by
$\cardw{A}$.
A morphism $f:(A,w) \to (B,v)$ is a function $f:A\to B$ which is
weight-preserving, i.e., a function $f$ such that $v(f(a)) = w(a)$
for all $a\in A$. By abuse of notation, we will often refer to
an object $(A,w)$ of $\catsets{R}$ simply as $A$, and write 
``$A\in \catsets{R}$.''

An \emph{$R$-weighted species}, 
\index{species!weighted}
or simply \emph{weighted species} if
$R$ is given, is a functor $F:\catfsets{}\to
\catsets{R}$. As before, $\sym{n}$ acts on $F[n]$, and we define
$\fixw F[\sigma]$ to be the sum of the weights of all elements
of $F[n]$ fixed by $\sigma$. This sum clearly depends only on
the cycle type of $\sigma$, so we can define $\fixw F[\lm]$
for any $\lm \vdash n$.
We define the \emph{weighted cycle
index series} of $F$ to be:

\begin{equation}
Z_F = \sum_{\lm} \fixw F[\lm] \pz{\lm}
\label{eq:wcycdef}
\end{equation}

Any species can be considered as a weighted species simply by
assigning $w(a)=1$ for all $a\in F[A]$, for any finite set
$A$. We then have $\fixw F[\lm] = \fix F[\lm]$ for any $\lm$,
and \eqref{eq:cycdef} and \eqref{eq:wcycdef} are identical. We will 
sometimes refer to a species, considered as a
weighted species in this fashion, as an \emph{ordinary} species.
\index{species!ordinary}

We will use the following notation, much of it standard in the
theory of species. By $\mathbf{0}$ (or simply $0$ when clear from the context)
we denote the \emph{empty species}, $\mathbf{0}[U]=\emptyset$ for all
finite sets $U$. By $\mathbf{1}$ (or $1$) we denote the \emph{empty set species},
\begin{equation*}
\mathbf{1}[U] = 
\begin{cases}   
\{U\} & \text{if } U=\emptyset\\
\emptyset & \text{otherwise}
\end{cases}
\end{equation*}
Similarly, we have the species $X$ of \emph{singletons},
\begin{equation*}
X [U] = 
\begin{cases}   
\{U\} & \text{if } |U|=1 \\
\emptyset & \text{otherwise}
\end{cases}
\end{equation*}
and the species $E_k$ of sets of cardinality $k$:
\begin{equation*}
E_k [U] = 
\begin{cases}   
\{U\} & \text{if } |U|=k \\
\emptyset & \text{otherwise}
\end{cases}
\end{equation*}
The \emph{uniform species} $E$, or \emph{species of sets},  is
given by
$E[U] = \{ U \}$ for all finite sets $U$. Thus,
$E = E_0 + E_1 + E_2 + \dots = 1 + X + E_2 + \dots$.

The combinatorial operations of sum ($+$), product ($\cdot$),
derivation ($'$), and substitution ($\circ$) on species induce corresponding
operations on their cycle index series---i.e., on symmetric functions.
The operations corresponding to $+$ and $\cdot$ are easily described (they
are the ordinary sum and product on symmetric functions), and the
operation corresponding to $'$ is $\partial\ / \partial p_1$.
That is,
$Z_{F'} = \partial Z_F / \partial p_1$,
where $Z_F$ is given in terms of the $p_i$ by \eqref{eq:cycdef}.
The operation corresponding to $\circ$ is well-known as well (see
\cite[Theorem 3]{spec1}); in our framework it corresponds to
the operation of plethysm of symmetric functions. (We explore substitution
and plethysm in more detail in Section~\ref{sec:comp-spec}.)

The Cartesian product 
\index{Cartesian product}
($\times$) of species $F$ and $G$ is defined by
$(F\times G)[U] = F[U]\times G[U]$ for all finite sets $U$. The corresponding
operation on symmetric functions is the internal, or Kronecker product. 
\index{Kronecker product}
Following \cite{lagrangeinv}, we will denote it by $\times$ as well. Thus,
\begin{equation*}
\sum_{\lm} a_{\lm} \pz{\lm} \times \sum_{\mu} b_{\mu} \pz{\mu}
= \sum_{\lm} a_{\lm} b_{\lm} \pz{\lm}
\end{equation*}
and for species $F$ and $G$, $Z_{F\times G} = Z_F \times Z_G$.

Closely related to the Cartesian product is the \emph{scalar product}
\index{scalar product}
$\langle \thinspace , \rangle$
(see \cite{lagrangeinv,anfunc}), given by 
$\langle F , G \rangle = (F\times G)|_{X=1}$. The
scalar product on species corresponds to the scalar product on
symmetric functions, which is given by:
\begin{equation}
\left\langle \sum_{\lm} a_{\lm} \pz{\lm}, \sum_{\mu} b_{\mu} \pz{\mu}
\right\rangle
=
\sum_{\lm} \frac{a_{\lm} b_{\lm}}{z_{\lm}}
\label{eq:scalar-prod}
\end{equation}

These notions generalize to species of several
sorts (see \cite{lagrangeinv} for details). For example, 
the cycle index of a 2-sorted 
species $F(X,Y)$ is given by:
\begin{equation}
Z_F = \sum_{\lm,\mu} \fix F[\lm,\mu] \pzm{x}{\lm} \pzm{y}{\mu}
\label{eq:2varcyc}
\end{equation}
(Here we follow the notation of~\cite{lagrangeinv} for symmetric functions
in several sets of variables: given $f\in \lam{}{}$, $f(y)$ denotes $f$ as
a symmetric function of the variables $y_1,y_2,\dots$, $f(z)$ denotes
$f$ in the variables $z_1,z_2,\dots$, etc.
Thus,
$p_1 (y) = y_1 + y_2 + y_3 + \dots$, $p_2 (z) = {z_1}^2 + {z_2}^2 +\dots$,
and so forth.)
For 2-sorted species
$F(X,Y)$, $G(X,Y)$, their
\emph{Cartesian product in $Y$} 
\index{Cartesian product!in $Y$}
is given by:
\begin{equation}
(F \times_Y G)[U,V]= \sum_{U_1 + U_2 = U} F[U_1,V] \times G[U_2,V]
\end{equation}
and the cycle index of $F \times_Y G$ is obtained by expressing
$Z_F$ and $Z_G$ in terms of the $p_{\lm} (y)$, with
coefficients in $\lamhat{x}$:

If
\begin{equation*}
\begin{aligned}
Z_F & = \sum_{\lm} a_{\lm} (x) \pzm{y}{\lm} \\
Z_G & = \sum_{\mu} b_{\mu} (x) \pzm{y}{\mu} 
\end{aligned}
\end{equation*}
for some  $a_{\lm} (x), b_{\mu} (x) \text{ in } \lamhat{x}$, then
\begin{equation*}
Z_{F \times_Y G} = \sum_{\lm} a_{\lm} (x) b_{\lm} (x) \pzm{y}{\lm}
\end{equation*}
The \emph{scalar product in $Y$}
\index{scalar product!in $Y$}
of $F$ and $G$ is then given by
$\langle F, G \rangle_{Y} = (F\times_Y G)|_{Y=1}$, and its cycle index
is obtained from $Z_{F \times_Y G}$ by setting each $p_{\lm} (y)$
equal to $1$. (Inspection of \eqref{eq:scalar-prod} shows that this
procedure is entirely analogous to that in the 1-sorted case.)

Finally, we recall that both the \emph{exponential generating function}
\index{generating function!exponential}
and the \emph{isomorphism-types generating function}
\index{generating function!isomorphism-types}
of a species
$F$ may be obtained from its cycle index series by certain
substitutions. In our framework, the exponential generating
function $F(x)$ is obtained by setting $p_1$ equal to $x$, and
$p_i$ equal to $0$ for $i>1$. The isomorphism-types generating function
$\widetilde{F}(x)$ is obtained by setting $p_i$ equal
to $x^i$ for all $i$. (For an $n$-sorted species, one performs these
substitutions in each variable. For example, the 
isomorphism-types generating function $\widetilde{F}(x,y)$
of a 2-sorted species
$F(X,Y)$ is obtained from the cycle index~\eqref{eq:2varcyc} by
setting $p_i (x)$ equal to $x^i$ and $p_i (y)$ equal to
$y^i$ for all $i$.)

\subsection{Burnside's Lemma}

For manipulating and calculating with cycle index series, a 
generalization of Burnside's
Lemma will be vital to us. We state it here.
(This result is not new; it is essentially a weighted
variation of Lemma 5 in~\cite{lagrangeinv}.)

Suppose a finite group $G\times H$ acts on a set $S\in \catsets{R}$,
and $w:S\to R$
is the weight function of $S$. By this we mean that $G\times H$
acts on the set $S$, and that for any $(g,h)\in G\times H$, the map
$s\mapsto (g,h)\cdot s$ is a morphism in $\catsets{R}$.
Let $\fixw{(g,h)}$ denote the sum of the weights
of all elements of $S$ fixed by $(g,h)$. The groups $G$ and $H$, 
considered as subgroups of $G\times H$, also act 
on $S$, and for $s\in S$, let $\orbit{G}{s}$ denote the orbit of
$s$ under the action of $G$, or simply the ``$G$-orbit'' of $s$.

\begin{lemma}
The weight function $w$ is constant on each $G$-orbit, allowing
us to define $w(\orbit{G}{s})$ to be $w(s)$. The group
$H$ acts on the set of $G$-orbits, and for any $h\in H$,
the sum of the weights of the $G$-orbits fixed by $h$ is:
\begin{equation}
\frac{1}{|G|} \sum_{g\in G} \fixw{(g,h)}
\end{equation}
\label{lem:orbwt}
\end{lemma}

\begin{proof}
Considering $G$ and $H$ as subgroups of $G\times H$ means identifying
$G$ with the set $G\times \{e_H\}$ and $H$ with the set $\{e_G\}\times H$,
where $e_G$ and $e_H$ denote the identity elements
of $G$ and $H$, respectively. Thus
the actions of $G$ and $H$ upon $S$ are given by
$g\cdot s = (g,e_H)\cdot s$,
$h\cdot s = (e_G,h)\cdot s$, for $g\in G$, $h\in H$, and $s\in S$. 
We note
that
$g\cdot (h\cdot s) = (gh)\cdot s = (hg) \cdot s = h\cdot (g\cdot s)$.
(By $gh$ we mean $(g,h)$, since
$(g,e_H)(e_G,h) = (g,h)$ in $G\times H$.)

It is immediate that $w$ is constant on each orbit under $G$, under
$H$, and under $G\times H$, by the definition of morphisms in
$\catsets{R}$.
To show that $H$ acts on
the set of $G$-orbits, define:
\begin{equation}
h\cdot \orbit{G}{s} = \orbit{G}{h\cdot s}
\label{eq:h-act}
\end{equation}

The action \eqref{eq:h-act} is well-defined, for if $\orbit{G}{s_1}
=\orbit{G}{s_2}$ then there exists $g\in G$ such that $ s_1
= g\cdot s_2$, and,
\begin{equation*}
\begin{aligned}
h\cdot \orbit{G}{s_1} & = h\cdot \orbit{G}{g\cdot s_2} \\
& = \orbit{G}{h\cdot (g\cdot s_2)} \\
& = \orbit{G}{(hg)\cdot s_2} \\
& = \orbit{G}{(gh)\cdot s_2} \\
& = \orbit{G}{g\cdot (h\cdot s_2)} \\
& = \orbit{G}{h\cdot s_2} \\
& = h\cdot \orbit{G}{s_2} 
\end{aligned}
\end{equation*}
That~\eqref{eq:h-act} \emph{is} a group action is clear: 
\begin{equation*}
\begin{aligned}
h_1 \cdot (h_2 \cdot
\orbit{G}{s} ) & = \orbit{G}{h_1\cdot (h_2\cdot s)} \\
& = \orbit{G}{(h_1 h_2) 
\cdot s} \\
 & = (h_1 h_2)\cdot \orbit{G}{s}
\end{aligned}
\end{equation*}

To calculate the sum of the weights of the $G$-orbits fixed by
$h\in H$, we consider the set $\mathcal{S}$ of
pairs $(g,s)$, where $g\in G$, $s\in S$, and $\inv{g}\cdot s = h\cdot s$.
Define the weight $w(g,s)$ of the pair $(g,s)$ to be $w(s)$.
Now $\inv{g}\cdot s = h\cdot s$ \myiff \ $(g,h)\cdot s = s$, i.e., 
\myiff \ $(g,h)$ fixes $s$. By the definition of $\fixw$, this means that for
a given $g\in G$, the sum of the weights of  pairs of the form $(g,s)$
$\mathcal{S}$
is simply $\fixw(g,h)$. Thus,
\begin{equation}
\sum_{(g,s)\in \mathcal{S}} w(g,s) = \sum_{g\in G} \fixw (g,h)
\label{eq:wsum1}
\end{equation}

We now calculate this sum of weights in a different way. For
a pair $(g,s)\in \mathcal{S}$, consider the number of choices for
$g$, given $s$. If $\orbit{G}{s}$ is not fixed by $h$, clearly there
are none (because if there exists $g\in G$ such that
$\inv{g}\cdot s = h\cdot s$, then $h\cdot s \in \orbit{G}{s}$,
and thus $h\cdot \orbit{G}{s} = \orbit{G}{s}$ by~\eqref{eq:h-act}).
If $h$ does fix $\orbit{G}{s}$ then $\orbit{G}{s} = \orbit{G}{h\cdot s}$,
so there exists $g_1\in G$ such that $g_1 \cdot s = h\cdot s$. Now,
$\inv{g}\cdot s = h\cdot s  \Leftrightarrow \inv{g}\cdot s = g_1 \cdot s 
\Leftrightarrow s = g g_1 \cdot s \Leftrightarrow g g_1 \in \stab{G}{s}$,
where $\stab{G}{s}= \{g\in G : g\cdot s = s \}$ denotes the $G$-stabilizer of $s$.
So the number of choices for $g$ is the number of $g$ such that
$g g_1 \in \stab{G}{s}$, which is simply $|\stab{G}{s}|$ (since
$g g_1 \in \stab{G}{s} \Leftrightarrow g\in \stab{G}{s} \inv{g_1}$,
and $|\stab{G}{s}\inv{g_1}|=|\stab{G}{s}|$). 

Let $\orbit{G}{s_1},\dots,\orbit{G}{s_k}$ be the $G$-orbits fixed by $h$.
Since there are $|\stab{G}{s}|$ choices for $g$ in the pair $(g,s)\in 
\mathcal{S}$, given $s$, we have:
\begin{equation*}
\begin{aligned}
\sum_{(g,s)\in \mathcal{S}} w(g,s) & = 
 \sum_{i=1}^{k} \sum_{s\in \orbit{G}{s_i}} |\stab{G}{s}| w(s) \\
& = \sum_{i=1}^{k} \sum_{s\in \orbit{G}{s_i}} \frac{|G|}{|\orbit{G}{s}|} w(s) \\
& = \sum_{i=1}^{k} |G| \sum_{s\in \orbit{G}{s_i}} \frac{w(s)}{|\orbit{G}{s}|} 
\end{aligned}
\end{equation*}
The weight function $w$ is constant on orbits, and $|\orbit{G}{s}| =
|\orbit{G}{s_i}|$ if $s$ is an element of $\orbit{G}{s_i}$, so,
\begin{equation}
\sum_{i=1}^{k} |G| \sum_{s\in \orbit{G}{s_i}} \frac{w(s)}{|\orbit{G}{s}|}  = 
|G| \sum_{i=1}^{k}  w(\orbit{G}{s_i})
\label{eq:wsum2}
\end{equation}
by the definition of the weight of an orbit.

Combining \eqref{eq:wsum1} and \eqref{eq:wsum2} completes the proof.
\end{proof}

\begin{corollary}[Burnside's Lemma, weighted form]
Suppose a finite group $G$ acts on a set $S$, and  $w$ is  a
weight function which is constant on orbits. Let $\fixw g$ denote
the sum of the weights of elements of $S$ fixed by $g\in G$.
Then the sum of the
weights of the orbits under the action of $G$ is
\begin{equation}
\frac{1}{|G|} \sum_{g\in G} \fixw g
\label{eq:burnside}
\end{equation}
\label{cor:burnside}
\end{corollary}

\begin{proof}
Apply Lemma~\ref{lem:orbwt} to the case where $H$ is a trivial group.
\end{proof}

Corollary~\ref{cor:burnside} is of course well-known 
(see for example (2.3.10)
in~\cite{graphenum}).

\subsection{Analytic functors}
\label{sec:anfunc}

In~\cite{anfunc}, Joyal associates to the species $F$ an \emph{analytic functor}
\index{analytic functor}
$F():\catsets{} \to \catsets{}$, where $\catsets{}$ is the category of sets and 
functions. His construction is
as follows. For a set $A$, consider the set $F[n] \times {A}^n$ of 
all pairs $(t,f)$, where $t\in F[n]$ and $f$ is a function from $[n]$ to $A$. 
The symmetric group $\sym{n}$ acts on these pairs via: 
\begin{equation}
\sigma \cdot (t,f) = (\sigma \cdot t,f \circ \inv{\sigma})
\label{eq:symact}
\end{equation}
where $\sigma \cdot t$ denotes
$F[\sigma](t)$. Let $F[n] \times {A}^n / \sym{n}$ be the set
of orbits under this action, and let $\orbit{\sym{n}}{t,f}$ be the orbit of the
pair $(t,f)$. The analytic functor $F()$ is then defined by:
\begin{equation}
F(A) = \sum_{n\geq 0} F[n]\times A^n / \sym{n}
\label{eq:functdef}
\end{equation}
for all sets $A$. (For a function $g:A\to B$, the map $F(g):F(A)\to F(B)$
is defined by $\orbit{\sym{n}}{t,f} \mapsto \orbit{\sym{n}}{t,g\circ f}$. It
is easily verified that $F(g)$ is well-defined and preserves composition.)

\begin{remark}
The analytic functor $F()$, written with parentheses, must not be 
confused with the species $F$, which is itself a functor and is often
written 
$F[\,]$.
\end{remark}

We generalize the notion of analytic functor to that of  
\emph{weighted analytic functor}
\index{analytic functor!weighted}
by modifying the category on which $F()$ is defined. For a weighted species
$F$, we construct the weighted analytic functor $F():\catsets{R}\to
\catsets{R}$ associated to a species $F$ as follows. For
$A\in \catsets{R}$, we have the set of pairs $(t,f)\in F[n]\times A^n$, and
the orbits $F[n]\times A^n / \sym{n}$, as 
above. Define the weight of $f$ to be:
\begin{equation}
w(f) = \prod_{a\in A} {w(a)}^{|\inv{f}(a)|}
= w(f(1)) w(f(2)) \dots
\label{eq:genweightdef}
\end{equation}
For any $\sigma\in \sym{n}$ and $a\in A$, the inverse images of
$a$ under $f$ and $f\circ \inv{\sigma}$ have the same size (since
$\inv{(f\circ \inv{\sigma})}(a) = \sigma (\inv{f}(a))$, and $\sigma$ is
a permutation), so the weight we have defined is constant on orbits
of $A^n$
under $\sym{n}$.  Define the weight of the pair $(t,f)$ to be
$w(t) w(f)$.
The
functor $F()$ is then defined by equation~\eqref{eq:functdef}, with
the weight $w(\orbit{\sym{n}}{t,f})$
of an orbit defined to be the weight of any of its elements.

\begin{definition}
For a weighted analytic functor $F$, and $A\in \catsets{R}$, we define the 
\emph{type-series}
of $F(A)$ to be:
\begin{equation}
Z_{F(A)} = \cardw{F(A)} = \sum_{\mathcal{O}\in F(A)} w(\mathcal{O})
\label{eq:type-series}
\end{equation}
\end{definition}

We now relate the type-series to the cycle index~\eqref{eq:wcycdef}. Consider
the case where $R= \mathbb{Q}[[ x_1,x_2,\dots]]$, the ring of formal
power series in the variables $x_1,x_2,\dots$, with coefficients in 
$\mathbb{Q}$. We consider the set
$\posints = \{1,2,\dots\}$ of positive integers to be an element of $\catsets{R}$
by assigning it the following weight function: for $i\in \posints$, let
$w(i) = x_i$. (For any $S\subseteq \posints$, we define a weight function
on $S$ by restriction.)
By~\eqref{eq:genweightdef}, the weight of a function $f:[n]\to
\posints$ is the following
monomial in the $x_i$:
\begin{equation}
w(f) = {x_1}^{|\inv{f}(1)|} {x_2}^{|\inv{f}(2)|} \dots
\label{eq:weightdef}
\end{equation}

We will often refer to $f$ as a ``coloring'' of the elements of $[n]$,
and think of $(t,f)$ as the $F$-structure $t$ together with a coloring of
the elements of its underlying set $[n]$.  An orbit under
the action $\sym{n}$---i.e., an element of $F(\posints)$---thus 
corresponds to an ``unlabeled, colored $F$-structure,''
or an ``isomorphism class of colorings of $F$-structures.''

In order to count the orbits under $\sym{n}$ by weight, we will use
Burnside's Lemma. To apply it, we must calculate the sum of the weights of 
pairs $(t,f)$ fixed by a permutation $\sigma
\in \sym{n}$. 

Suppose $\sigma$ has cycle type $\lm = (1^{m_1} 2^{m_2} \dots)$. If
$\sigma$ fixes a pair $(t,f)$ then by \eqref{eq:symact}, $\sigma$
must fix $t$ (that is, $F[\sigma](t) = t$), and $f$ must be constant on
each cycle of $\sigma$. If a cycle $c$ of length $l$ receives color $i$
(that is, if $f(j) = i$ for all $j$ in the cycle $c$), this results
in a factor of ${x_i}^l$ in the weight of $(t,f)$. The colors can
be assigned arbitrarily, and there are $m_l$ cycles of length $l$, the
sum of the weights of the colorings of a given $F$-structure $t$
fixed by $\sigma$ is:
\begin{equation}
(x_1 + x_2 + \dots)^{m_1} ({x_1}^2 + {x_2}^2 + \dots)^{m_2} 
({x_1}^3 + {x_2}^3 + \dots)^{m_3} \dots
\label{eq:colorsum}
\end{equation}
since each term in the expansion of this product corresponds to
the weight of a specific coloring.

The product \eqref{eq:colorsum} is immediately expressible in
terms of power sum symmetric functions as ${p_1}^{m_1} {p_2}^{m_2} \dots
= p_{\lm}$. Thus, we find that the sum of the weights of pairs
$(t,f)$ fixed by $\sigma$ of cycle type $\lm$ is:
\begin{equation*}
\fixw \sigma = \fixw F[\sigma] p_{\lm}
\end{equation*}
Applying Corollary~\ref{cor:burnside},
the sum of the weights of orbits under the action \eqref{eq:symact}
is:
\begin{equation}
\frac{1}{n!} \sum_{\sigma \in \sym{n}} \fixw F[\sigma] p_{\lm}
\label{eq:symsum1}
\end{equation}
where $\lm$ is the cycle type of $\sigma$. Using the fact that
there are $n! / z_{\lm}$ permutations of cycle type $\lm$ in $\sym{n}$, 
we can rewrite
\eqref{eq:symsum1} as:
\begin{equation}
\sum_{\lm \vdash n} \fixw F[\lm] \pz{\lm}
\label{eq:symsum2}
\end{equation}

In light of \eqref{eq:symsum2},
the combinatorial interpretation of \eqref{eq:cycdef} is now clear:
the cycle index of a species $F$, considered as a formal series in 
the $x_i$, counts
by weight all isomorphism classes of colorings of $F$-structures. (The
sum \eqref{eq:symsum2} counts only those on underlying sets of cardinality
$n$.) The connection with the type-series is:
\begin{equation}
Z_{F(\posints)} = Z_F
\label{eq:cyc-equal}
\end{equation}
where the $Z_F$ on the right-hand side is the 
weighted cycle index~\eqref{eq:wcycdef},
considered as a formal series in the $x_i$.

\begin{remark}
The connection between \eqref{eq:cycdef} and \eqref{eq:symsum1}, with
the $p_i$ replaced by formal variables, is of course well-known (see for example
\cite[Proposition 13]{spec1}). In fact, the equality \eqref{eq:cyc-equal}
could also have been proven by P{\'o}lya theory. We recall that given a
subgroup $G$ of $\sym{n}$, its cycle index polynomial $Z(G)$ is given by
(using the $p_i$ in place of formal variables):
\begin{equation*}
Z(G)=\frac{1}{|G|} \sum_{g\in G} p_{\lm}
\end{equation*}
where $\lm$ is the cycle type of $g$.  By Proposition 13 of
\cite{spec1}, the cycle index $Z_F$ of a species $F$ can be expressed as:
\begin{equation}
Z_F = \sum_{t} Z(\aut{t})
\label{eq:aut-sum}
\end{equation}
where the sum is over a
set of representatives of all isomorphism classes of $F$-structures, and
$\aut{t}$ denotes the automorphism group of an $F$-structure $t$.
By P{\'o}lya's theorem (see, for example, (2.4.16) in
\cite{graphenum} for a statement), 
\index{P{\'o}lya's Theorem}
replacing each $p_i$ by
${x_1}^i + {x_2}^i + \dots + {x_n}^i$ in \eqref{eq:aut-sum} gives a formal
series which counts by weight the inequivalent colorings of
$F$-structures, colored in $n$ colors. Since this holds for any positive
$n$, we obtain the equality \eqref{eq:cyc-equal}.
\end{remark}

\subsection{Composition of species}
\label{sec:comp-spec}

As an application of our interpretation of the cycle index,
we
provide a proof of the well-known expression for the cycle index of the
composition of two species. We first recall the operation $\circ$ 
of \emph{plethysm}
(see~\cite[Section I.8]{symfunc}) on symmetric functions. Plethysm
satisfies $(fg)\circ h = (f\circ h)(g\circ h)$, $(f+g)\circ h =
f\circ h + g\circ h$ for all $f,g,h \in \lam{}{}$. The plethysm
$p_n \circ g$ is obtained from the expression for $g$ (in
terms of the $x_i$) by replacing each $x_i$ with ${x_i}^n$. For example, it 
follows immediately that $p_n \circ p_m = p_{nm}$ for all $n,m \in \mathbb{N}$.
(We will deal with plethysm in more detail, in the context of \lamrings,
in Section~\ref{sec:inplethsym}.)

\begin{lemma}
Let $F$ and $G$ be ordinary species, considered as weighted
species, with associated weighted analytic functors
$F()$, $G()$. Then:
\begin{equation}
Z_{F(G(\posints))} = Z_{F(\posints)} \circ Z_{G(\posints)}
\end{equation}
\label{lem:compfunc}
\end{lemma}

\begin{proof}
By definition, an element of $F(G(\posints))$ is the orbit under
$\sym{n}$, for some $n$, of a pair $(t,f)$, where $t\in F[n]$ and
$f:[n]\to G(\posints)$. The weight of the pair $(t,f)$ is, 
by~\eqref{eq:genweightdef}:
\begin{equation*}
w(t,f) = 
\prod_{\mathcal{O}\in G(\posints)} {w(\mathcal{O})}^{|\inv{f}(\mathcal{O})|}
\end{equation*}

Suppose $\sigma\in \sym{n}$ fixes $(t,f)$, where $\sigma$ has cycle type
$\lm = (1^{m_1} 2^{m_2} \dots)$. Then $\sigma$ must fix $t$, and since
$f\circ \inv{\sigma} = f$, $f$ must be constant on the cycles of $\sigma$.
Suppose the value of $f$ on a cycle $c$ of length $l$ in $\sigma$ is
$\mathcal{O} \in G(\posints)$. The contribution to the weight of $(t,f)$
will be a factor of $w(\mathcal{O})^l$. Since $w(\mathcal{O})$ is
a monomial (by~\eqref{eq:weightdef}), $w(\mathcal{O})^l$ can be
obtained from $w(\mathcal{O})$ by replacing each $x_j$ in
$w(\mathcal{O})$ with ${x_j}^l$. Thus, given $t$, the sum of the weights
of pairs $(t,f)$ fixed by $\sigma$ is:
\begin{equation}
(p_1 \circ Z_{G(\posints)})^{m_1} (p_2 \circ Z_{G(\posints)})^{m_2}\dots
= p_{\lm} \circ Z_{G(\posints)}
\label{last13}
\end{equation}
since each term in the expansion of this product corresponds to
the weight of such a pair. Therefore, $\fixw \sigma$ is simply
$\fix F[\sigma]$ times the product~\eqref{last13}, and by~\eqref{eq:burnside},
the sum of the weights of the orbits under $\sym{n}$ is:
\begin{equation*}
\frac{1}{n!} \sum_{\sigma \in \sym{n}} p_{\lm} \circ Z_{G(\posints)}
\end{equation*}
where $\lm$ is the cycle type of $\sigma$.

Summing over all $n$ and applying the fact that there are $n!/z_{\lm}$
permutations of cycle type $\lm$ in $\sym{n}$, we find that:
\begin{equation*}
\begin{aligned}
Z_{F(G(\posints))} & = \sum_{\lm} \frac{\fix F[\lm]}{z_{\lm}} p_{\lm} \circ
 Z_{G(\posints)} \\
& = Z_{F(\posints)} \circ Z_{G(\posints)}
\end{aligned}
\end{equation*}
which completes the proof.
\end{proof}

\begin{lemma}
Let $F$ and $G$ be weighted species with associated weighted analytic functors
$F()$, $G()$, and let $(F\circ G)()$ be the weighted analytic functor
associated to the species $F\circ G$. There is a natural isomorphism between
the functors $(F\circ G)()$ and $F(G())$.
\label{lem:compfunc-spec}
\end{lemma}

\begin{proof}
We must construct
a weight-preserving bijection between $(F\circ G)(A)$ and
$F(G(A))$ for any finite set $A$. Consider first $(F\circ G)(A)$.

An element of $(F\circ G)(A)$ consists of the orbit of
a pair $(t,f)$, where
$t\in (F\circ G)[m]$ and $f:[m]\to A$, for some $m$. By the definition
of the composition of two species, $t$ consists of:
\begin{myenumerate}
\item A partition $U_1 + \dots + U_n = [m]$, for some $n$
\item $G$-structures $\alpha_1,\dots, \alpha_n$, with $\alpha_i \in
 G[U_i]$
\item An $F$-structure $\alpha \in F[\{\alpha_1,\dots, \alpha_n\}]$
\end{myenumerate}

We can thus map the orbit of $(t,f)$ to the orbit of an element 
$(t',f')$ as follows.
Let $\gamma : \{\alpha_1,\dots, \alpha_n\} \to [n]$ be any bijection,
and let $\overline{\sigma}\in \sym{n}$ be the permutation defined by
$\overline{\sigma}(i) = k \Leftrightarrow \gamma(\alpha_k) = i$.
For $1\leq i \leq n$, let $m_i = |U_i|$, and
let $\gamma_i: U_i \to [m_i]$ also be
bijections. Define $t' = F[\gamma](\alpha)$, and for $1\leq i \leq n$,
let $f': [n] \to G(A)$ be defined by 
$f'(i) = \orbit{\sym{m_k}}{G[\gamma_k](\alpha_k), f\circ \inv{\gamma_k}}$
for $1\leq i \leq n$, where $k=\overline{\sigma}(i)$.

Let $\Gamma$ be the map we have constructed. We first show that
$\Gamma$ is well-defined. 

The orbit of $(t',f')$ does not depend on the choice of the bijections
$\gamma$ and $\gamma_1,\dots,\gamma_n$, for any other bijections will
be of the form $\sigma \circ \gamma$, $\sigma_1 \circ \gamma_1,
\dots, \sigma_n \circ \gamma_n$ for some permutations $\sigma \in
\sym{n}$, $\sigma_i \in \sym{m_i}$, $1\leq i \leq n$. Replacing $\gamma$
by $\sigma \circ \gamma$ means replacing $\overline{\sigma}$ with
$\overline{\sigma}\circ \inv{\sigma}$, since
$(\overline{\sigma}\circ \inv{\sigma})(i) = k \Leftrightarrow
\gamma(\alpha_k) = \inv{\sigma}(i) \Leftrightarrow
(\sigma \circ \gamma)(\alpha_k) = i$. 
Thus,
using these
new bijections maps the orbit of $(t,f)$ to the orbit of
$({t_1}',{f_1}')$, where ${t_1}' = F[\sigma \circ \gamma](\alpha)$
and ${f_1}' (i) = 
\orbit{\sym{m_i}}{G[\sigma_k \circ \gamma_k](\alpha_k), f\circ \inv{\gamma_k}
\circ \inv{\sigma_k}}$,
where $k = (\overline{\sigma}\circ \inv{\sigma})(i)$.
But by the definition~\eqref{eq:symact} of the actions
of the various symmetric groups, it is immediate that
$\orbit{\sym{n}}{t',f'} = \orbit{\sym{n}}{{t_1}',{f_1}'}$.

The orbit of $(t',f')$ also does not depend on the choice of
the representative $(t,f)$, for any other representative will be of
the form $(\sigma \cdot t, f\circ \inv{\sigma})$ for some
$\sigma \in \sym{m}$. The element $\sigma\cdot t$ will consist of:
\begin{myenumerate}
\item The partition $\sigma (U_1) + \dots + \sigma (U_n) = [m]$
\item The $G$-structures $G[\restrict{\sigma}{U_1}](\alpha_1),\dots, 
 G[\restrict{\sigma}{U_n}](\alpha_n)$
\item The $F$-structure $F[\delta](\alpha)$, where
$\delta$ is the map $\alpha_i \mapsto 
G[\restrict{\sigma}{U_i}](\alpha_i)$ for $1\leq i \leq n$
\end{myenumerate}

Let $\gamma '$ be the bijection $\gamma \circ \inv{\delta}$,
and ${\gamma_i}'$ be $\gamma_i \circ \inv{(\restrict{\sigma}{U_i})}$ for
$1\leq i \leq n$. Using these bijections, one verifies immediately
that $\Gamma$ maps the orbit of  $(\sigma \cdot t, f\circ \inv{\sigma})$ to
that of $(t',f')$.

To show that $\Gamma$ is injective, suppose the orbit of $(u,g)$ also maps to the
orbit of $(t',f')$. Then $u$ must consist of:
\begin{myenumerate}
\item A partition $\myprime{U_1} + \dots + \myprime{U_n} = [m]$
\item $G$-structures $\myprime{\alpha_1},\dots, \myprime{\alpha_n}$, 
 with $\myprime{\alpha_i} \in G[\myprime{U_i}]$.
\item An $F$-structure $\myprime{\alpha} \in 
 F[\{\myprime{\alpha_1},\dots, \myprime{\alpha_n}\}]$
\end{myenumerate}

Let $\myprime{\gamma}, \myprime{\gamma_1},\dots,\myprime{\gamma_n}$ be the
bijections used in the definition of $\Gamma$ to calculate the image of
$(u,g)$, and $\myprime{\overline{\sigma}}$ the permutation
corresponding to $\myprime{\gamma}$---i.e., 
$\myprime{\overline{\sigma}}(i) = k \Leftrightarrow 
\myprime{\gamma}(\myprime{\alpha_k}) = i$. Without loss of generality---since
we have shown that the bijections used in the definition of $\Gamma$ can
be chosen arbitrarily---we may assume that $\myprime{\overline{\sigma}}$
is the identity. Thus, $\myprime{\gamma}(\myprime{\alpha_i}) = i$ for
$1\leq i \leq n$. Let $\myprime{m_i} = |\myprime{U_i}|$.

Using the bijections 
$\myprime{\gamma}, \myprime{\gamma_1},\dots,\myprime{\gamma_n}$
in the definition of $\Gamma$ gives a representative $(u',g')$ of the orbit
which $(u,g)$ maps to, and since this is also the orbit to which
$(t,f)$ maps, this representative must be of the form $(\sigma \cdot t',
f' \circ \inv{\sigma})$ for some $\sigma \in \sym{n}$. Therefore,
$F[\myprime{\gamma}](\myprime{\alpha}) = F[\sigma \circ \gamma](\alpha)$.
Again because the choice of bijections is arbitrary, we may replace $\gamma$
by $\sigma \circ \gamma$ and assume 
$F[\myprime{\gamma}](\myprime{\alpha}) = F[\gamma](\alpha)$.
Recalling that $\myprime{\gamma}(\myprime{\alpha_i}) = i$ and
that $\gamma(\alpha_k) = i \Leftrightarrow \overline{\sigma}(i)=k$,
we see that:
\begin{equation*}
\begin{aligned}
\inv{\gamma} \circ \myprime{\gamma}: 
\myprime{\alpha_i} & \mapsto \alpha_{\overline{\sigma}(i)} \\
F[\inv{\gamma} \circ \myprime{\gamma}](\myprime{\alpha}) &= \alpha
\end{aligned}
\end{equation*}

Moreover, because $(u',g') = (t',f')$, we have that $g' = f'$, and thus:
\begin{equation*}
\begin{aligned}
\orbit{\sym{ \myprime{m_i}}}{ G[\myprime{\gamma_i}](\myprime{\alpha_i}),
 g\circ \inv{(\myprime{\gamma_i})} }  = &
\orbit{\sym{ m_k}}{ G[\gamma_k](\alpha_k),
 f\circ \inv{\gamma_k} } \\
 & \text{where } k=\overline{\sigma}(i) 
\end{aligned}
\end{equation*}
The $\gamma_k$ may be chosen arbitrarily up to the action of a permutation,
so we may in fact assume that $G[\myprime{\gamma_i}](\myprime{\alpha_i})
= G[\gamma_k](\alpha_k)$ and
$g\circ \inv{(\myprime{\gamma_i})} = f\circ \inv{\gamma_k}$
for $1\leq i \leq n$.

Consider the maps $\inv{\gamma_k} \circ
\myprime{\gamma_i}$ from $\myprime{U_i}$ to $U_k$. Glueing these maps together
yields a permutation $\tau$ of $[m]$, and it is straightforward to verify
that $\tau \cdot (u,g) = (t,f)$.

To show that $\Gamma$
is surjective we construct, given $(t',f')$,
an element $(t,f)$ which maps to it. For $1\leq i \leq n$,
choose $\myprime{\alpha_i}$ in the orbit $f'(i)$. Thus, $\myprime{\alpha_i}$ 
consists of a pair $(\myprime{t_i},\myprime{f_i})$, with $\myprime{t_i}\in G[m_i]$,
$\myprime{f_i}:[m_i]\to A$, for some $m_i$. Let $U_1+ \dots + U_n$
be any partition of $m=m_1+\dots m_n$ such that $|U_i|=m_i$, and
choose $\gamma_i:U_i \to [m_i]$ to be any bijections.
Let $\alpha_i = G[\inv{\gamma_i}](\myprime{\alpha_i})$,
and define $f_i = \myprime{f_i}\circ \gamma_i$. Let $\gamma :
\{\alpha_1,\dots,\alpha_n \}$ be any bijection, and let
$t=F[\inv{\gamma}](t')$. Define $f:[m]\to A$ by glueing
together the $f_i$. Then by construction, the orbit of
$(t,f)$ maps to the orbit of $(t',f')$ under $\Gamma$.

The verification that $\Gamma$ is weight-preserving, and of
naturality, is straightforward.
\end{proof}

\begin{theorem}
Let $F$ and $G$ be ordinary species. Then
\begin{equation}
Z_{F\circ G} = Z_F \circ Z_G
\label{eq:comp-spec}
\end{equation}
\end{theorem}

\begin{proof}

This follows immediately from~\eqref{eq:cyc-equal} and 
Lemmas~\ref{lem:compfunc} and
\ref{lem:compfunc-spec}, since:
\begin{equation*}
\begin{aligned}
Z_{F\circ G} &= Z_{(F\circ G)(\posints)}\\
 &= Z_{F(G(\posints))} \\
 &= Z_{F(\posints)} \circ Z_{G(\posints)} \\
 &= Z_F \circ Z_G
\end{aligned}
\end{equation*}
\end{proof}

\section{The inner plethysm of symmetric functions}
\label{sec:inplethsym}

We first review some facts about representation theory. By a
\emph{representation} we shall always mean a complex,
finite-dimensional representation, i.e., a homomorphism from
a finite group $G$ to $\genlin (n,\cmplx)$ for some $n\geq 0$.
 We take our notation for \lamrings\ from \cite{lamrings1}. In particular,
given a \lamring\ $R$, $\psi^{n}:R\to R$ will denote the $n$th Adams
operation of $R$.

Given a finite group $G$, $\repring (G)$ will denote the \emph{representation
ring} of $G$, that is, the integer span of all isomorphism classes
of representations of $G$, with addition
given by direct sum, and multiplication by tensor product. For a 
representation $\rho$ of $G$, $[\rho] \in \repring (G)$ will denote its
isomorphism class.
The representation
ring is a \lamring \ via the operations:

\begin{equation*}
\lambda^{n}([\rho]) = [\wedge^{n} (\rho)]
\end{equation*}
where $\wedge^{n}$ denotes the $n$th exterior power.

We  
denote by $\centfunc{G}{}$ the ring of all central functions of $G$, that
is, the ring of all functions from $G$ to $\cmplx$ which are constant on
conjugacy classes. (Addition and multiplication are defined via pointwise
addition and multiplication in $\cmplx$.)
The \lamring \ structure of $\centfunc{G}{}$ is given
by its Adams operations,
\begin{equation}
(\adams{k}(c))(g) = c(g^k) 
\label{linkg}
\end{equation}
for all $c\in \centfunc{G}{}$ and $g\in G$.

There is a map $\chi:\repring (G) \to \centfunc{G}{}$ which
sends  $[\rho] \in \repring (G)$ to the character of $\rho$, $\chi_{\rho}
= \trace \circ \rho$. (Here $\trace$ denotes trace.) As shown in \cite{lamrings1},
$\chi$ sends
$\repring (G)$ isomorphically  onto its image, which we will
denote here by $\charring{G}$, the \emph{character ring} of $G$. The irreducible
characters (that is, images of irreducible representations) form a 
$\mathbb{Z}$-basis for $\charring{G}$, and a $\cmplx$-basis for $\centfunc{G}{}$.

If $G$ is the symmetric group $\sym{n}$, we can say more. As shown
in \cite{symfunc}, $\charring{\sym{n}}$ is isomorphic to $\Lambda^{n}$, the ring
of symmetric functions homogeneous of degree $n$ (where the
multiplication in $\Lambda^{n}$ is taken to be the Kronecker product). 
The isomorphism is
given by the characteristic map $\ch$,

\begin{equation}
\ch(\tau) = \sum_{\lambda \vdash n}
\tau_{\lambda} \pz{\lambda}
\label{linkd}
\end{equation}
where $\tau_{\lambda}$ denotes the value of $\tau \in \charring{\sym{n}}$ on 
permutations
of cycle type $\lambda$.

Given a \lamring \ $R$, the ring $\lam{}{}$ of symmetric function
acts upon $R$ in a natural way. The construction is given 
by Knutson (see~\cite[page 46]{lamrings1}) as follows. For $f\in \lam{}{}$,
we can express $f$ as a polynomial $F$ in the elementary symmetric
functions $e_i$: $f=F(e_1,e_2,\dots)$. The symmetric function
$f$ then becomes a natural operation on $R$ via
$r\mapsto F(\lm (r), {\lm}^2 (r), {\lm}^3 (r), \dots)$ for all $r\in R$.

\begin{definition}
For $f$ and $r$ as above, we will denote 
$F(\lm (r), {\lm}^2 (r), {\lm}^3 (r), \dots)$ by $f[r]$.
\label{linke}
\end{definition}

\begin{remark}
We can define an action of $\lamhat{}$ upon a \lamring\ similarly.
However, an element $f\in \lamhat{}$ is in general a
formal series $F(e_1,e_2,\dots)$, and thus
$f[r]$ may not be defined.
\end{remark}

\begin{lemma}
Let $a$ and $b$ be elements of a \lamring \ $R$, and
$n,m\geq 1$. Then,
\begin{myenumerate}
\item $ p_n [r] = \adams{n}(r) $
\item $ p_n [a+b] = p_n [a] + p_n [b]$
\item $ p_n [ab] = p_n [a] p_n [b]$
\item $ p_n [p_m [a]] = p_{nm} [a]$
\end{myenumerate}
\end{lemma}

\begin{proof}
Item (1) follows immediately from the definition of the Adams
operations (see~\cite[page 47]{lamrings1}). Items (2) -- (4) then
follow from the properties of $\adams{n}$ (see~\cite[page 48]{lamrings1}).
\end{proof}

\begin{example}
The most natural application of Definition~\ref{linke} is the case in which
$R$ is $\lam{}{}$  itself (since $\lam{}{}$ is the free \lamring \ on one
generator). In this case, for $f,g \in \Lambda$, $f[g]$ is simply
the plethysm $f \circ g$.
\end{example}

We describe the inner plethysm of symmetric functions via Definition~\ref{linke}
 as well.
Set $R=\charring{\sym{n}}$  and denote by $\inpleth{}$ the action
of $\lam{}{}$  on $\lam{n}{}$  induced by identifying $\charring{\sym{n}}$
 and $\lam{n}{}$
via the isomorphism \eqref{linkd}. For $f \in \Lambda$, $g\in \Lambda^{n}$,
$f\inpleth{} g$ is the \emph{inner plethysm}
\index{inner plethysm!of symmetric functions}
 of $f$ and $g$. To extend
$\inpleth{}$ to an action $\Lambda \times \Lambda \to \Lambda$, we define
$f \inpleth{} (g_1+g_2) = f\inpleth{} g_1 + f\inpleth{} g_2$ for
$g_1\in \Lambda^{n}, g_2 \in \Lambda^{m}$ with $n \neq m$.

Inner plethysm corresponds to the composition of a representation of
$\sym{n}$ with a representation of a general linear group (indeed, this is
how it was originally defined by Littlewood~\cite{littlewood}).

Let $\phi: \genlin (m,\cmplx) \to \genlin (n,\cmplx)$ be a representation.
As in \cite{GLn}, if $\phi$ is a polynomial representation, we associate to 
$\phi$ its character $f_{\phi}$,
that is, the symmetric function $f(x_1,\dots,x_n)$ such that:
\begin{equation*}
f(\theta_1,\dots,\theta_n)=\trace (\phi (A))
\end{equation*}
for any $A\in GL(m,\cmplx)$,
where $\theta_1,\dots,\theta_n$ are the eigenvalues of $A$.

We note that $f_{\phi}$ is a symmetric function 
in a finite number of variables. However, by abuse of notation, we
will also denote by $f_{\phi}$ any element
$f\in \Lambda$ such that 
$f(x_1,\dots,x_n,0,0,\dots) = f_{\phi}(x_1,\dots,x_n)$.

\begin{theorem}
Let $\rho:\sym{n} \to \genlin (m,\cmplx)$ and $\phi: \genlin (m,\cmplx) \to
 \genlin (q,\cmplx)$
be representations. Then,
\begin{equation*}
\chi_{\phi \circ \rho} = f_{\phi}[\chi_{\rho}]
\end{equation*}
in the \lamring \ $\charring{\sym{n}}$, the action of $\lam{}{}$ 
on $\charring{\sym{n}}$
being given by Definition~\ref{linke}.
\end{theorem}

\begin{proof}
Fix $\sigma \in \sym{n}$, and let $M= \rho (\sigma)$.
Suppose $M$ has eigenvalues $\theta_1,\dots,\theta_m$.
Then $\chi_{\rho} (\sigma) = \theta_1 + \dots + \theta_m$, by
definition, and since the trace of the $k$th
power of a matrix is the sum of the $k$th powers of the
eigenvalues, $\chi_{\rho} (\sigma^k) = 
\theta_1{}^{k} + \dots + \theta_m{}^{k}$.
By \eqref{linkg} and Definition~\ref{linke}, this means:
\begin{equation*}
\begin{aligned}
(p_k [\chi_{\rho}]) (\sigma) & = \theta_1{}^{k} + \dots + \theta_m{}^{k} \\
& = p_k (\theta_1,\dots,\theta_k)
\end{aligned}
\end{equation*}

It follows immediately from Definition~\ref{linke} 
that $(p_{\lambda} [\chi_{\rho}]) (\sigma) =
p_{\lambda} (\theta_1,\dots,\theta_k)$ for any partition $\lambda$.

Now express $f_{\phi}$ in terms of the $p_{\lambda}$: $f_{\phi} =
\sum_{\lambda} a_{\lambda} p_{\lambda}$ for some coefficients
$a_{\lambda} \in \mathbb{Q}$. Then,
\begin{equation}
\begin{aligned}
\trace \phi (M) & = \sum_{\lambda} a_{\lambda} p_{\lambda} 
(\theta_1,\dots,\theta_k) \\
& = \sum_{\lambda} a_{\lambda} (p_{\lambda}[\chi_{\rho}]) (\sigma) \\
& = f_{\phi} [\chi_{\rho}] (\sigma)
\end{aligned}
\label{linkf}
\end{equation}
This proves the theorem. 
\end{proof}

It follows immediately from our definition of inner plethysm that
with $\rho$ and $\phi$ as in \eqref{linkf},
\begin{equation*}
\ch (\chi_{\phi \circ \rho}) = f_{\phi} \inpleth{} \ch (\chi_{\rho})
\end{equation*}
Thus, the inner plethysm of symmetric functions
corresponds to the composition of a representation of a symmetric
group with a representation of a general linear group.

\section{The inner plethysm of species}

We will now define an operation of inner plethysm on species,
such that the cycle index of the inner plethysm of two species
is the inner plethysm of their cycle indices. In effect, this will
give a combinatorial interpretation to the operation of inner plethysm
on symmetric functions. We give two constructions of the operation,
one using combinatorial operations on species, the other using
the analytic functors introduced in Section~\ref{sec:anfunc}.

\subsection{Combinatorial construction}
\label{sec:inplethspec}

In order to define the inner plethysm of species, we first study a certain
map $\Phi$ from 1-sorted species to 2-sorted species.
Let $\mathcal{H}(X,Y)$ be the 2-sorted
species defined by letting $\mathcal{H}[U,V]$ be the set of all functions from
$U$ to $V$. For a 1-sorted species $F$, define $\Phi(F)(X,Y)$ by:
\begin{equation}
\Phi(F) [U,V] = \mathcal{H}[U,F[V]]
\label{eq:phidef}
\end{equation}

\begin{remark}
We note that $\Phi(0)=E(Y)$. This is because $\Phi(0)[U,V]=\emptyset$
if $U \neq \emptyset$, but if $U=\emptyset$ and $V$ is any set then
$\Phi(0)[U,V]$ contains a single element: the empty function
from $\emptyset$ to $\emptyset$.
\end{remark}

\begin{lemma}
\begin{equation*}
\fix \Phi (F)[\beta,\sigma]
= \prod_{k \geq 1} (\fix F[\sigma^k])^{\beta_k}
\end{equation*}
where $\beta_k$ denotes the number of cycles of length $k$ of
the permutation $\beta$.
\label{lem:Phi-fix}
\end{lemma}

\begin{proof}
Suppose a function $f$ is fixed by $\Phi (F)[\beta,\sigma]$. Then
$f$ satisfies $F[\sigma] \circ f \circ \inv{\beta} = f$. The image
of every element of a cycle $c$ of $\beta$ of length $l$ under $f$ is 
thus determined
by the image of some element $i\in c$; moreover, $F[\sigma^l]$ must
fix $f(i)$. So the total number of possibilities for $f$ is:
\begin{equation*}
\prod_{c} \fix F[\sigma^{l(c)}]
\end{equation*}
where the product is over all cycles $c$ of $\beta$, and
$l(c)$ is the length of $c$. This is equivalent to the statement of
the lemma.
\end{proof}

The cycle index series of $\Phi (F)$ is thus:
\begin{equation}
\begin{aligned}
Z_{\Phi (F)} &= \sum_{\lambda, \mu} \biggl( \prod_{k \geq 1} 
{\fix F[\sigma^k]}^{\beta_k} \biggr) \frac{p_{\lambda}(x)}{z_{\lambda}}
\frac{p_{\mu}(y)}{z_{\mu}} \\
&= \sum_{\mu} \exp \biggl( \sum_{i \geq 1} \frac{\fix F[\sigma^i] p_{i}(x)}{i}
\biggr) \frac{p_{\mu}(y)}{z_{\mu}}
\end{aligned}
\label{eq:Phicyc}
\end{equation}
where $\beta_k$ is the number of $k$-cycles of a permutation of
cycle type $\lambda$, and $\sigma$ is a permutation of cycle type $\mu$. 

\begin{remark}
Lemma~\ref{lem:Phi-fix} could also have been proven by noting that
$\mathcal{H}(X,Y) = E(E(X)\cdot Y)$ and using this to calculate
$Z_{\mathcal{H}}$, then applying the techniques in
\cite{functcomp} (which easily extend to the multi-variable case) to
calculate $Z_{\Phi}$.
\end{remark}

We now define the inner plethysm of 1-sorted species, which we will
denote~$\inpleth{}$.

\begin{definition}
For 1-sorted species $F$, $G$,
\begin{equation*}
(F \inpleth{} G)(Y) = \langle F(X), \Phi (G) (X,Y) \rangle_X
\end{equation*}
\label{def:inplethdef}
\end{definition}
\index{inner plethysm!of species}

We note that $F \inpleth{} G$ is defined \myiff\ $F$ is
\emph{strictly finite}
\index{species!strictly finite}
(in the terminology of \cite{virtspec}), i.e.,
\myiff\ there exists an integer $N$ such that $F[n] = \emptyset$
whenever $n>N$.

\begin{theorem}
\begin{equation}
Z_{F \inpleth{} G} = Z_F \inpleth{} Z_G
\label{eq:pleth-spec}
\end{equation}
\end{theorem}

\begin{proof}
By linearity, it is sufficient to consider the case where $G$ is
homogeneous of degree $m$, for some $m$. Then $G$ gives a permutation 
representation
of $\sym{m}$ with character $\chi$,
\begin{equation*}
\chi (\sigma) = \fix G[\sigma]
\end{equation*}
for $\sigma \in \sym{m}$.

Suppose $\lambda$ is a partition, $\lambda = (1^{\beta_1} 2^{\beta_2} \dots)$.
Then, in the \lamring \ $\charring{\sym{m}}$, with the action of
$\Lambda$ on $\charring{\sym{m}}$ given by Definition~\ref{linke},
\begin{equation*}
(p_{\lambda}[\chi])(\sigma) = 
\chi (\sigma)^{\beta_1} \chi (\sigma^2)^{\beta_2} \chi (\sigma^3)^{\beta_3}
\dots
\end{equation*}
for $\sigma \in \sym{m}$.
Therefore, by \eqref{eq:Phicyc},
\begin{equation*}
Z_{\Phi(G)} = \sum_{\lambda} \sum_{\mu \vdash m}
(p_{\lambda}[\chi])(\sigma)
\pzm{x}{\lambda} \pzm{y}{\mu}
\end{equation*}
where $\sigma$ has cycle type $\mu$.

Let 
\begin{equation*}
Z_F = \sum_{\lambda} f_{\lambda} \pzm{x}{\lambda}
\end{equation*}
Then
\begin{equation}
Z_F \times_{X} Z_{\Phi(G)} = 
\sum_{\lambda} 
 \biggl(\sum_{\mu \vdash m} (p_{\lambda}[\chi])(\sigma)
 \pzm{y}{\mu} \biggr)
  f_{\lambda} \pzm{x}{\lambda}
\label{last7}
\end{equation}
where $\sigma$ has cycle type $\mu$.

Setting $p_{\lambda}(x)=1$ in \eqref{last7} for all $\lambda$ shows that
the species
$F\inpleth{} G$ gives a representation of $\sym{m}$ whose character
$\tau$ is:
\begin{equation}
\begin{aligned}
\tau (\sigma) & = \sum_{\lambda} 
f_{\lambda} \frac{(p_{\lambda}[\chi])(\sigma)}{z_{\lambda}} \\
& = (Z_F [\chi])(\sigma)
\end{aligned}
\label{last8}
\end{equation}
for $\sigma \in \sym{m}$.
But $\chi$ is the character of the representation of
$\sym{m}$ corresponding to $G$, i.e. $\ch(\chi) = Z_G$. So by the
definition of inner plethysm, \eqref{last8} means that
$\ch (\tau) = Z_F \inpleth{} Z_G$.
Since $\tau$ is the character of the representation of $\sym{m}$
corresponding to $F\inpleth{} G$, $\ch (\tau) = Z_{F\inpleth{} G}$, and
the theorem is proved. 
\end{proof}

\subsection{Construction by analytic functors}

Given species $F$ and $G$, we have associated
functors $F()$ and $G()$ such that the functorial composition of
the latter corresponds to the composition, in the species-theoretic sense,
of the former. Since $F$ and $G$ are themselves functors, we can consider
other compositions;  $F(G[\,])$, for example. 
We note that $F(G[\,])$ is a  species \myiff\ $F$ is strictly finite; in
this case, $F(U)$ is finite for any finite set $U$, and thus
$F(G[A])$ is finite for any finite set $A$.

\begin{lemma}
There is a natural isomorphism between the functors $F(G[\,])$ and
$F\inpleth{} G$, where $\inpleth{}$ is the operation on species
constructed in Section~\ref{sec:inplethspec}.
\end{lemma}

\begin{proof}
Let $A$ be a finite set. An element of $F(G[A])$ consists, for some $n$,
 of the
orbit under $\sym{n}$ of a pair $(t,f)$, where $t\in F[n]$,
and $f:[n] \to G[A]$.  The action of $\sym{n}$ is
$\sigma \cdot (t,f) = (\sigma \cdot t, f\circ \inv{\sigma})$.

By Definition~\ref{def:inplethdef}, an element of $(F\inpleth{} G)[A]$
consists, for some $n$, of the orbit under $\sym{n}$ of a
pair $(t,f)$, where
$t\in F[n]$, and $f\in \Phi(G)[n,A]$. By~\eqref{eq:phidef},
this means $f$ is a function from $[n]$ to $A$. The action of $\sym{n}$ is
$\sigma \cdot (t,f) = (\sigma \cdot t, f\circ \inv{\sigma})$.

The lemma follows immediately.
\end{proof}

Thus, we find that
the operation of inner plethysm arises extremely naturally in the context of
analytic functors. Indeed, focusing on analytic functors as the 
objects of study, we could define the inner plethysm
of $F$ and $G$ to be $F(G[\,])$. To show that this point of view is 
useful, we must show that this definition would provide a means of
calculating the cycle index of $F\inpleth{} G$.

\begin{theorem}
\begin{equation}
Z_{F \inpleth{} G} = Z_F \inpleth{} Z_G
\end{equation}
\end{theorem}

\begin{proof}

Let $H$ be the species $F(G[\,])$. By~\eqref{eq:cyc-equal}, to calculate
the cycle index of $F(G[\,])$, it will be sufficient to calculate the sum
of the weights of the elements of $H(\posints)$.

An element of $H(\posints)$ consists of the orbit under
$\sym{n}$ of a pair $(t,f)$,
where $t\in H[n]$, $f:[n]\to \posints$, for some $n$.
This means $t\in F(G[n])$, and so $t$ is the orbit under $\sym{m}$
of $(t_1,f_1)$,
where $t_1\in F[m]$, $f_1:[m]\to G[n]$, for some $m$. 

We therefore have an action of $\sym{m} \times \sym{n}$ on the
set $\mathcal{S}$ of
triples $(t_1,f_1,f)$, given by:
\begin{equation}
(\tau,\sigma) \cdot (t_1,f_1,f) =
 (\tau \cdot t_1, G[\sigma]\circ f_1 \circ \inv{\tau}, f\circ \inv{\sigma})
\label{last14}
\end{equation}
for $(\tau,\sigma)\in \sym{m} \times \sym{n}$. The symmetric groups
$\sym{m}$ and $\sym{n}$ act separately as subgroups of
$\sym{m} \times \sym{n}$, and an orbit under $\sym{m}$ corresponds
to a pair $(t,f)$. The symmetric group $\sym{n}$ acts on these
orbits, and an orbit of $\sym{m}$-orbits under $\sym{n}$
corresponds to an element of $H(\posints)$. 

To calculate the
sum of the weights of orbits of pairs $(t,f)$ under $\sym{n}$ by
Burnside's Lemma, we must calculate $\fixw \sigma$ for $\sigma\in
\sym{n}$. By the previous paragraph, this amounts to calculating the
sum of the weights of $\sym{m}$-orbits of $\mathcal{S}$ fixed
by $\sigma$. By Lemma~\ref{lem:orbwt}, we have:
\begin{equation}
\fixw \sigma = \frac1{m!} \sum_{\tau\in \sym{m}} \fixw (\tau,\sigma)
\end{equation}
where $\fixw (\tau,\sigma)$ is the sum of the weights of elements
of $\mathcal{S}$ fixed by $(\tau,\sigma)$.

By \eqref{last14}, if $(t_1,f_1,f)$ is fixed by $(\tau,\sigma)$, we
must have:
\begin{myenumerate}
\item $F[\tau](t_1) = t_1$
\item $G[\sigma]\circ f_1 \circ \inv{\tau} = f_1$
\item $f\circ \inv{\sigma} = f$
\end{myenumerate}
So the number of choices for $t_1$ in such a triple is simply
$\fix F[\tau]$, by item (1) above. By item (3), the weight
of the triple is $p_{\lm}$, where the partition $\lm$ is the cycle
type of $\sigma$---by the same argument used to 
derive~\eqref{eq:colorsum}. 

We examine item (2).
If $f_1$ satisfies $G[\sigma]\circ f_1 \circ \inv{\tau} = f_1$, and $c$
is a cycle of $\tau$ of length $l$, with $i\in c$, it is immediate that 
the image of any
element of $c$ under $f_1$ is determined by the image $f_1 (i)$.
Moreover,  $G[{\sigma}^l]$ must fix $f_1 (i)$.
So we see that the total number of choices for $f_1$ is:
\begin{equation}
\prod_{c} \fix G[{\sigma}^{l(c)}]
\label{last15}
\end{equation}
where the product is over all cycles $c$ of $\tau$, and $l(c)$ is
the length of $c$. Let $\chi \in \charring{\sym{n}}$ be the
element $\chi (\sigma) = \fix G[\sigma]$. Then~\eqref{last15}
is clearly equal to $p_{\mu} [\chi](\sigma)$, 
where $\mu$ is the cycle type of $\tau$ and the action of
$\lam{}{}$ on $\charring{\sym{n}}$ is given by Definition~\ref{linke}.

So we see that 
$\fixw (\tau,\sigma) = \fix F[\tau] p_{\mu} [\chi](\sigma) p_{\lm}$,
and thus that for
$\sigma\in \sym{n}$, the sum of the weights of
$\sym{m}$-orbits fixed by $\sigma$ is:
\begin{equation}
\begin{aligned}
\fixw \sigma &= \frac{1}{m!} \sum_{\tau \in \sym{m}} 
 \fix F[\tau] p_{\mu}[\chi](\sigma) p_{\lm} \\
 & \quad \text{where } \mu,\,\lm \text{ are the cycle types of }
  \tau,\,\sigma \text{ respectively,} \\
 & = \sum_{\mu \vdash m} \frac{\fix F[\mu] p_{\mu}[\chi](\sigma)}{z_{\mu}}
     p_{\lm}
\end{aligned}
\end{equation}
Summing over all $m$, we find that for the
action of $\sigma$ on $H(\posints)$,
$\fixw \sigma = Z_F [\chi](\sigma) p_{\lm}$, and so the sum
of the weights of the orbits under $\sym{n}$ is:
\begin{equation}
\frac{1}{n!} \sum_{\sigma\in \sym{n}} Z_F [\chi](\sigma) p_{\lm}
 = 
 \sum_{\lm \vdash n} Z_F [\chi](\sigma) \pz{\lm}
\end{equation}
where on the right-hand side, $\sigma$ is any permutation of cycle
type $\lm$. Summing over all $n$ gives the cycle index of
$F\inpleth{} G$, and by the definition of the inner plethysm of
symmetric functions, we find that it is precisely 
$Z_F \inpleth{} Z_G$.
\end{proof}

\section{The inner plethysm in \textit{Y}}
\label{sec:inplethY}

We recall that a 2-sorted species $F$ is a functor
$F: {\catfsets{}}^2 \to \catfsets{}$, with associated cycle index series
$Z_F$:
\begin{equation}
Z_F = \sum_{\lm,\mu} \fix F[\lm,\mu] \pzm{x}{\lm} \pzm{y}{\mu}
\end{equation}

In direct analogy with Section~\ref{sec:anfunc}, we define an 
\emph{$R$-weighted, $2$-sorted species} to be a functor 
$F: {\catfsets{}}^2 \to \catsets{R}$, and associate to it
a weighted analytic functor $F(): {\catsets{R}}^2 \to
\catsets{R}$:
\begin{equation}
F(A,B) = \sum_{m,n \geq 0} F[m,n]\times A^m \times B^n / {\sym{m}\times \sym{n}}
\end{equation}
where the action of $\sym{n}\times \sym{m}$ is
$(\tau,\sigma) \cdot
(t,f,g) = ((\tau,\sigma)\cdot t, f\circ \inv{\tau},g\circ \inv{\sigma})$.
We define the weight $w(t,f,g)$ to be $w(t) w(f) w(g)$, and
associate to $F(A,B)$ its \emph{type-series},
\begin{equation}
Z_{F(A,B)} = \sum_{\mathcal{O}\in F(A,B)} w(\mathcal{O})
\label{eq:2vartype}
\end{equation}

The relation between the  cycle index and  type-series
is entirely analogous to that in the one-variable case. Let 
$R=\mathbb{Q}[[x_1,x_2,\dots,y_1,y_2,\dots]]$, let $\posints_{x}$
denote $\posints$ with with weight function $w(i)=x_i$, and let
$\posints_{y}$ denote $\posints$ with the weight function $w(i) = y_i$.
Then,
\begin{equation}
Z_{F(\posints_{x},\posints_{y})} = Z_F
\label{eq:2varcyc-equal}
\end{equation}
where  $Z_F$ on the right-hand side is the cycle index~\eqref{eq:2varcyc},
considered as a formal series in the $x_i$ and $y_i$.
The proof is entirely analogous to that of \eqref{eq:cyc-equal}.

An $n$-sorted species $F$, its cycle index, and its type-series, are defined
in the obvious way, and the analog of~\eqref{eq:2varcyc-equal} is easily
verified.

So we see that, as in the one-variable case, we can consider the cycle
index either as a symmetric function in the underlying variables, or
as a series in the $p_i (x)$, $p_i (y)$, etc. We can also, however, consider
it as a series in the $p_i (y)$ with coefficients in
$\mathbb{Q}[[x_1, x_2, \dots]]$. This interpretation will prove useful.

\subsection{Combinatorial description}
\label{sec:inplethY-description}

In Section~\ref{sec:inplethY-construction} we will construct an operation
$\inpleth{Y}$, the \emph{inner plethysm in $Y$}, which we describe here.
Intuitively, $F\inpleth{Y} G$ is similar to the
substitution $F(G(X,Y))$. Whereas an $F(G(X,Y))$-structure consists
of an $F$-enriched set of $G$-structures, however, an 
$F\inpleth{Y} G$-structure can be said to consist of ``an
$F$-enriched set of $G$-structures which share the same $Y$s.'' We
illustrate with an example.

\begin{figure}

\begin{center}
\epsfig{file=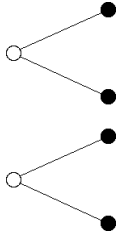,height=1.197in}
\end{center}

\caption{The isomorphism class of $G$-structures}
\label{fig:G-struct}
\end{figure}

Consider the case where $G = E_2 (X \cdot E_2 (Y))$. An isomorphism class
(in fact, there is only one) of $G$-structures is shown in 
Figure~\ref{fig:G-struct}. The white points are considered to be of
sort $X$, the black of sort $Y$.

\begin{figure}

\begin{center}
\epsfig{file=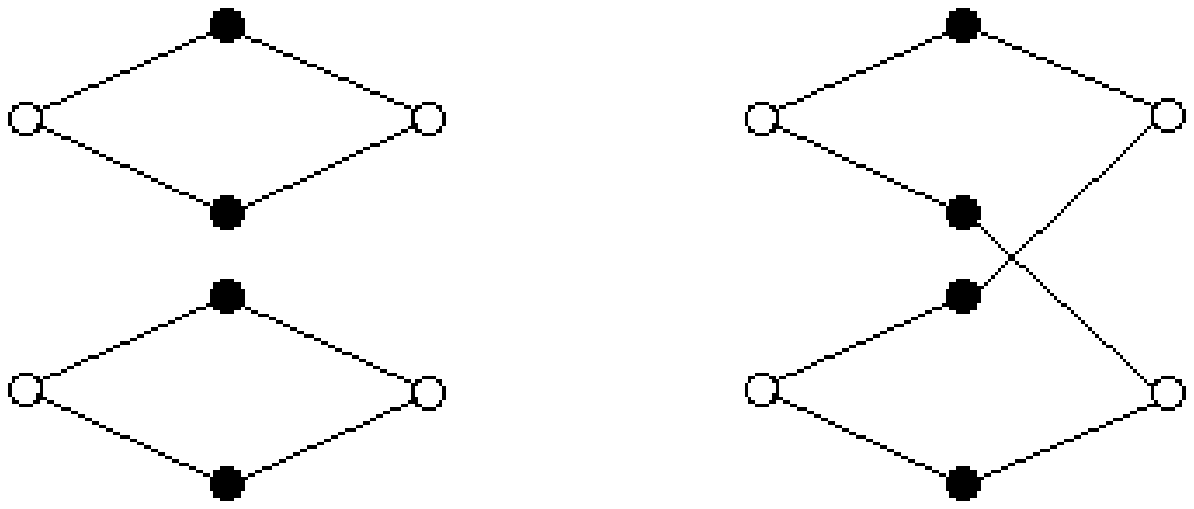,height=0.988in}
\end{center}
\caption{The isomorphism classes of $F\inpleth{Y}G$-structures}
\label{fig:FG-struct}
\end{figure}

Let $F=E_2$. We would expect there to be 2 isomorphism classes of
$F\inpleth{Y}G$-structures (see Figure~\ref{fig:FG-struct}), since
an $F\inpleth{Y}G$-structure consists of a set of 2 $G$-structures which
share the same $Y$-points. That this is, in fact, the case, we verify
once we have set up the necessary machinery.

\subsection{Construction}
\label{sec:inplethY-construction}

The standard notation for a 2-sorted species with underlying points of
sorts $X$ and $Y$ is $F(X,Y)$. Since this is easily confused with our
notation for the analytic functor associated to $F$, we will denote
the analytic functor by $F_{XY}$, and make the following definition.

\begin{definition}
For a 2-sorted species $F(X,Y)$, we define
\begin{equation}
\begin{gathered}
F_{X[Y]}  :   \catsets{R} \times \catfsets{}   \to \catsets{R} \\
F_{X[Y]} (A,B)   = \sum_{n\geq 0} F[n,B] \times A^n / \sym{n}
\end{gathered}
\label{eq:functdefX}
\end{equation}
\end{definition}

We will refer to $F_{X[Y]}$ as the \emph{analytic functor in $X$}
associated to $F$, or simply as a \emph{partial analytic functor}.
\index{analytic functor!partial}
(The action of $\sym{n}$ on $F[n,B] \times A^n$ in
\eqref{eq:functdefX} is the obvious one: $\sigma\cdot (t,f)
= ( F[\sigma,1_{B}](t), f\circ \inv{\sigma})$, where
$1_B$ is the identity function on $B$.)

For any set $A\in \catsets{R}$,
the analytic functor in $X$ gives a functor 
$F_{X[Y]}(A,-):\catfsets{}\to \catsets{R}$, i.e., a weighted species,
the underlying points of which are of sort $Y$. If $A=\posints_{x}$,
the cycle index of this species
is a symmetric function in $y_1,y_2,\dots$, with coefficients in
$\mathbb{Q}[[x_1,x_2,\dots]]$. Considered as a symmetric function in
$x_1,x_2,\dots,y_1,y_2,\dots$, it is precisely the cycle 
index~\eqref{eq:2varcyc}, as the following slightly more general result
will show.

\begin{lemma}
For the weighted, 1-sorted species 
$F_{X[Y]}(\posints_{x},-):\catfsets{}\to \catsets{R}$, we have:
\begin{equation}
\fixw F_{X[Y]}(\posints_{x},\sigma) =
\sum_{\lm} \fix F[\lm, \mu] \pzm{x}{\lm}
\end{equation}
where the sum is over all partitions $\lm$, and $\mu$ is the cycle
type of $\sigma$.
\label{lem:fixF}
\end{lemma}

\begin{proof}
Let $\sigma$ be an element of $\sym{n}$. By \eqref{eq:functdefX}, an
element of $F_{X[Y]}(\posints_{x},n)$ is the orbit under
$\sym{m}$, for some $m$, of a pair $(t,f)$, with $t\in F[m,n]$, 
$f:[m]\to \posints_{x}$. The weight of the pair is 
$w(f) = f(1) f(2)\dots$.
Thus, we have an action of $\sym{m}\times \sym{n}$ on pairs $(t,f)$:
\begin{equation*}
(\tau,\sigma)\cdot (t,f) = (F[\tau,\sigma](t),f\circ \inv{\tau})
\end{equation*}
and $\fixw F_{X[Y]}(\posints_{x},\sigma)$ is  the sum of the
weights of $\sym{m}$-orbits fixed by $\sigma$. Now, by the
argument used to derive~\eqref{eq:colorsum}, the sum of the weights
of functions $f$ such that
$f\circ \inv{\tau} = f$ is  $p_{\lm}(x)$, where $\lm$ is the cycle
type of $\tau$. So by Lemma~\ref{lem:orbwt}, the sum
of the weights of $\sym{m}$-orbits fixed by $\sigma$ is:
\begin{equation*}
\frac{1}{m!} \sum_{\tau\in \sym{m}} \fix F[\tau,\sigma] p_{\lm}(x)
 = \sum_{\lm \vdash m} \fix F[\lm,\mu] \pzm{x}{\lm}
\end{equation*}
where $\lm$ is the cycle type of $\tau$, and $\mu$ is the cycle type
of $\sigma$. 

Summing over all $m$ completes the proof.
\end{proof}

Suppose now that $F$ is a 1-sorted species, and $G(X,Y)$ is a
2-sorted species. We consider the composition
$F(G_{X[Y]} (-,-)):\catsets{R}\times \catfsets{}\to \catsets{R}$.

\begin{lemma}
With $F$ and $G$ as above, there exists a 2-sorted species $H(X,Y)$ such
that $H_{X[Y]} = F(G_{X[Y]} (-,-))$.
\label{lem:plethexist}
\end{lemma}

\begin{proof}
We  construct $H$ as follows. For finite sets $U$, $V$, we define
an element of $H[U,V]$ to be:
\begin{myenumerate}
\item A partition $U_1 + \dots + U_n = U$
\item $G$-structures $\alpha_1,\dots,\alpha_n$, with $\alpha_i \in G[U_i,V]$
\item An $F$-structure $\alpha \in F[\{\alpha_1,\dots,\alpha_n\}]$
\end{myenumerate}

For $A\in \catsets{R}$, $B\in \catfsets{}$,
an element of $H_{X[Y]}(A,B)$ consists of the orbit under $\sym{n}$,
for some $n$, of a pair $(t,f)$, with $t\in H[n,B]$, and $f:[n]\to A$.

Now 
an element of $F(G_{X[Y]} (A,B))$
consists of the orbit under $\sym{n}$, for some $n$, of a pair
$(t',f')$, with $t'\in F[n]$, $f':[n]\to G_{X[Y]} (A,B)$. Thus,
for $i\in [n]$, $f'(i)$ consists of the orbit under $\sym{m_i}$, for
some $m_i$, of a pair $({t_i}',{f_i}')$, where ${t_i}' \in G[m_i,B]$,
${f_i}' : [m_i]\to A$.

We can thus map the orbit of $(t,f)$ to the orbit of an element 
$(t',f')$ as follows.
Let $\gamma : \{\alpha_1,\dots, \alpha_n\} \to [n]$ be any bijection,
and let $\overline{\sigma}\in \sym{n}$ be the permutation defined by
$\overline{\sigma}(i) = k \Leftrightarrow \gamma(\alpha_k) = i$.
For $1\leq i \leq n$, let $m_i = |U_i|$, and
let $\gamma_i: U_i \to [m_i]$ also be
bijections. Define $t' = F[\gamma](\alpha)$, and for $1\leq i \leq n$,
let $f': [n] \to G_{X[Y]} (A,B)$ be defined by 
$f'(i) = \orbit{\sym{m_k}}{G[\gamma_k,1_{B}](\alpha_k), f\circ \inv{\gamma_k}}$
for $1\leq i \leq n$, where $k=\overline{\sigma}(i)$, and
$1_B$ is the identity function on $B$.

The verification that this gives a natural isomorphism between functors
is straightforward, and entirely analogous to the proof 
of Lemma~\ref{lem:compfunc-spec}.
\end{proof}

\begin{definition}
For species $F$, $G$ as above, we define the \emph{inner plethysm in $Y$}
\index{inner plethysm!in $Y$}
of $F$ and $G$ to be the species $H$ constructed in the proof of
Lemma~\ref{lem:plethexist}, and denote it $F \inpleth{Y} G$.
\end{definition}

In order to calculate the cycle index of $F \inpleth{Y} G$, we first
give a slight generalization of the \lamring \ of central functions of
a group.

\begin{definition}
For a group $G$ and a \lamring \ $R$, we define the ring 
$\centfunc{G}{R}$ of central
functions from $G$ to $R$ to be the set of all functions from $G$
to $R$ which are constant on conjugacy classes. The Adams operations
on $\centfunc{G}{R}$ are given by:
\begin{equation}
(\psi^k (f))(g) = \psi_{R}^k (f(g^k))
\end{equation}
for all $f\in \centfunc{G}{R}$, $g\in G$. Here $\psi_{R}^k$
denotes the $k$th Adams operation of the \lamring\ $R$.
\label{def:centfunc-lamring}
\end{definition}

By Lemma~\ref{lem:fixF}, in order to calculate the desired cycle
index, it will be sufficient to calculate
$\fixw (F\inpleth{Y} G)(\posints_{x},\sigma)$ for a permutation
$\sigma$.

\begin{lemma}
For species $F$ and $G$ as above, define 
$\chi: \sym{n}\to \lam{}{x}$ by:
\begin{equation*}
\chi (\sigma) = \fixw G_{X[Y]} (\posints_{x},\sigma)
= \sum_{\lm} \fix G[\lm, \mu] \pzm{x}{\lm}
\end{equation*}
where $\mu$ is the cycle type of $\sigma\in \sym{n}$. Then:
\begin{equation}
\fixw (F\inpleth{Y} G)(\posints_{x},\sigma) =
 Z_F [\chi](\sigma)
\end{equation}
where the action of $\lam{}{}$ on $\centfunc{\sym{n}}{\lam{}{x}}$
is given by Definition~\ref{linke}.
\label{lem:inplethYcalc}
\end{lemma}

\begin{proof}
By definition, an element of $(F\inpleth{Y} G)(\posints_{x},n)$
is the orbit under $\sym{m}$, for some $m$, of a pair
$(t,f)$, with $t\in F[m]$, $f:[m]\to G_{X[Y]} (\posints_{x},n)$.
The symmetric group $\sym{n}$ acts on $G_{X[Y]} (\posints_{x},n)$,
and so we have an action of $\sym{m}\times \sym{n}$ on the set
of pairs $(t,f)$. We wish to calculate 
$\fixw (F\inpleth{Y} G)(\posints_{x},\sigma)$, the sum of the weights
of $\sym{m}$-orbits fixed by $\sigma$.

Suppose a pair $(\tau,\sigma)\in \sym{m}\times \sym{n}$ fixes a pair
$(t,f)$. Then $\tau$ must fix $t$ (i.e., $F[\tau](t) = t$), and
$f$ must satisfy
\begin{equation}
G_{X[Y]}(\posints_{x},\sigma) \circ f \circ \inv{\tau} = f
\label{last17}
\end{equation}
Consider a cycle $c$ of $\tau$ of length $l$, with $i\in c$.
The image of any element of $c$ under $f$ is determined by the image
of $i$, and all such images have the same weight, so the cycle $c$
contributes a factor of $w(f(i))^l$ to the weight of $f$.

Since $c$ has length $l$, \eqref{last17} implies that
$G_{X[Y]}(\posints_{x},\sigma^l)$ must fix $f(i)$. Thus,
$\fixw G_{X[Y]}(\posints_{x},\sigma^l)$
counts by weight all possible choices for $f(i)$.
Since $f(i)$ is an element of $G_{X[Y]} (\posints_{x},n)$, 
$w(f(i))$ is a monomial
in the variables $x_1, x_2, \dots$, and $w(f(i))^l$
is obtained from $w(f(i))$ by replacing each variable $x_j$ with
${x_j}^l$.

Replacing each $x_j$ with ${x_j}^l$ in a symmetric
function $g$ amounts to calculating the plethysm $p_l \circ g$. So
we see that the sum of the weights of functions $f$
which satisfy~\eqref{last17} is:
\begin{equation*}
\prod_{c} p_{l(c)} \circ \fixw G_{X[Y]}(\posints_{x},\sigma^{l(c)})
= p_{\lm} [\chi] (\sigma)
\end{equation*}
where the product is over all cycles $c$ of $\tau$, and $\lm$ is the
cycle type $\tau$.

Thus, by Lemma~\ref{lem:orbwt}, the sum of the weights of $\sym{m}$-orbits 
fixed by $\sigma$ is:
\begin{equation*}
\begin{aligned}
\frac{1}{m!}\sum_{\tau \in \sym{m}} \fixw (\tau,\sigma)
 & = \frac{1}{m!}\sum_{\tau \in \sym{m}} \fix F[\tau] p_{\lm} [\chi] (\sigma) \\
 & = \sum_{\lm \vdash m} \fix F[\lm] \frac{p_{\lm} [\chi] (\sigma)}{z_{\lm}}
\end{aligned}
\end{equation*}

Summing over all $m$ completes the proof.
\end{proof}

Lemma~\ref{lem:inplethYcalc} provides a method for calculating the
cycle index of $F \inpleth{Y} G$, which can be described succinctly by
defining an operation $\inpleth{Y}: \lam{}{} \times \lam{}{xy} \to
\lam{}{xy}$ as follows. Suppose $f\in \lam{}{}$, and
$g\in \lam{}{xy}$ is homogeneous of degree $n$ in $y$, that is,
\begin{equation*}
g= \sum_{\mu\vdash n} a_{\mu}(x) \pzm{y}{\mu}
\end{equation*}
where $a_{\mu}\in \lam{}{x}$. Define $\chi: \sym{m} \to \lam{}{x}$
by $\chi (\sigma) = a_{\mu}(x)$, where $\mu$ is the cycle type of 
$\sigma$. Then $\chi \in \centfunc{\sym{n}}{\lam{}{x}}$, and we can
define:
\begin{equation*}
f\inpleth{Y} g= \sum_{\mu\vdash n} f[\chi] (\sigma) \pzm{y}{\mu}
\end{equation*}
where $\sigma$ has cycle type $\mu$, and the action of $\lam{}{}$
on $\centfunc{\sym{n}}{\lam{}{x}}$ is given by Definition~\ref{linke}.
(The \lamring\ structure of $\centfunc{\sym{n}}{\lam{}{x}}$ is given
by Definition~\ref{def:centfunc-lamring}.)
We extend $\inpleth{Y}$ to a mapping $\lam{}{} \times \lam{}{xy} \to
\lam{}{xy}$ by linearity. Lemma~\ref{lem:inplethYcalc} can then be
expressed as:
\begin{equation}
Z_{F \inpleth{Y} G} = Z_F \inpleth{Y} Z_G
\label{eq:inplethY-spec}
\end{equation}

We define operations $\inpleth{X}$, $\inpleth{XY}$, $\inpleth{XZ}$, and so
forth, in the obvious way. To define $\inpleth{XZ}$, for example, we define
$G_{[X]Y[Z]}$ for a 3-sorted species $G$ by:
\begin{equation*}
G_{[X]Y[Z]} (A,B,C) = \sum_{n\geq 0} G[A,n,C] \times B^n / \sym{n}
\end{equation*}
and for a 1-sorted species $F$, we define $F\inpleth{XZ} G$
to be the species $H$ such that
$H_{[X]Y[Z]} = F(G_{[X]Y[Z]}(-,-,-))$. The analogs of
Lemmas~\ref{lem:fixF}, \ref{lem:plethexist}, and \ref{lem:inplethYcalc},
are easily verified.

We now return to the example given in Section~\ref{sec:inplethY-description}.
 We have 
$Z_F = h_2 = \frac{1}{2}({p_1}^2 + p_2)$, and,
\begin{equation*}
\begin{aligned}
Z_G  =\:  & h_2 \circ (p_1 (x) h_2 (y)) \\
     = \:  & \frac{1}{8} \pp{1}{y}{4} \pp{1}{x}{2} + 
      \frac{1}{4} \pp{1}{y}{2} \pp{1}{x}{2} \pp{2}{y}{} +
      \frac{1}{8} \pp{1}{x}{2} \pp{2}{y}{2} + \\
      & \frac{1}{4} \pp{2}{y}{2} \pp{2}{x}{} +
      \frac{1}{4} \pp{4}{y}{} \pp{2}{x}{} 
\end{aligned}
\end{equation*}
Thus, by Lemma~\ref{lem:inplethYcalc},
\begin{equation*}
\begin{aligned}
Z_F \inpleth{Y} Z_G  = \: &
\frac{1}{4} \pp{4}{y}{} \pp{2}{x}{2} + \frac{1}{4} \pp{4}{y}{} \pp{4}{x}{}
  + \frac{3}{8} \pp{2}{y}{} \pp{1}{y}{2} \pp{2}{x}{2} + \\
 &  \frac{1}{8} \pp{2}{y}{} \pp{1}{y}{2} \pp{1}{x}{4} +
  \frac{7}{16} \pp{2}{y}{2} \pp{2}{x}{2} + 
   \frac{1}{16} \pp{2}{y}{2} \pp{1}{x}{4} + \\
 &   \frac{1}{4} \pp{2}{y}{2} \pp{1}{x}{2} \pp{2}{x}{} +
  \frac{1}{16} \pp{1}{y}{4} \pp{2}{x}{2}
 + \frac{3}{16} \pp{1}{y}{4} \pp{1}{x}{4}
\end{aligned}
\end{equation*}
(This calculation was performed with the aid of \maple.)

Setting $p_i (x) = x^i$, $p_i (y) = y^i$ for all $i$ in
$Z_F \inpleth{Y} Z_G$
yields $2 x^4 y^4$, indicating that there are 2 unlabeled
$F\inpleth{Y}G$-structures on  4 points of sort $X$ and 4 points of
sort $Y$, as expected.

\section{Polynomial maps}

We denote by $\spec{n}$ the set of all $n$-sorted
species, and by $\mspec{n}$ the set of all $n$-sorted 
molecular species. We consider two species $A$ and $B$ to be equal,
and write ``$A=B$,'' if they are naturally isomorphic. The set $\spec{n}$ is
a half-ring, and not a ring, since it has operations $+$ and $\cdot$, but
no operation of subtraction. Just as the integers can be constructed from the
natural numbers, the ring
$\vspec{n}$  of all $n$-sorted \emph{virtual species}
\index{species!virtual}
can be constructed
from $\spec{n}$. As shown in \cite{virtspec}, we
have the following half-ring and ring isomorphisms:

\begin{align*}
\spec{n}   & \cong   \mathbb{N}[[\mspec{n}]] \\
\vspec{n}  & \cong   \mathbb{Z}[[\mspec{n}]]
\end{align*}

Similarly, by $\specf{n}$ $\vspecf{n}$ we denote the sets of
strictly finite species and virtual species, respectively. We have 
(again, see \cite{virtspec}):

\begin{align*}
\specf{n}   & \cong   \mathbb{N}[\mspec{n}] \\
\vspecf{n}  & \cong   \mathbb{Z}[\mspec{n}]
\end{align*}

The operations $+$ and $\cdot$ extend from $\spec{n}$ to $\vspec{n}$
by construction. In this section we develop techniques---essentially
those used in \cite{virtspec}---to
extend others, such as $\circ$,
$\inpleth{}$, $\inpleth{Y}$, and to prove that the corresponding 
identities for symmetric functions (\eqref{eq:comp-spec},
\eqref{eq:pleth-spec}, and \eqref{eq:inplethY-spec}, in the
case of these three) remain valid for virtual species.

\subsection{Polynomial maps on species}
\label{sec:polymaps-spec}

We will use the notation $(a_i)_{i\in I}$ to denote an indexed collection of
(not-necessarily-distinct) objects $a_i$ (indexed by the elements of 
a some set $I$). Strictly speaking, $(a_i)_{i\in I}$ is shorthand for
the function defined on $I$ which has the value $a_i$ at $i$
for every $i\in I$. Thus, $(a_i)_{i \in \mathbb{N}}$, with 
$a_i \in \mathbb{R}$,
denotes a sequence of real numbers, and
if $R$ is an $n$-sorted species with
molecular series $\sum_{N \in \mspec{n}} a_N N$, then
$(a_N)_{N\in \mspec{n}}$ is the collection of coefficients of $R$.

Given a species $R$, $[R]$ will denote the collection of coefficients
of the molecular series of $R$.

We will often be concerned with functions from species to species.
A function $f:\spec{n} \to \spec{m}$ determines a collection of functions 
$(f_M)_{M\in \mspec{m}}$ via its action on molecular series:

\begin{equation}
f(R) = \sum_{M\in \mspec{m}}
f_M ([R]) M
\label{linka}
\end{equation}

\begin{definition}
If the functions $f_{M}$ in \eqref{linka} are polynomials in the 
coefficients $[R]=(a_N)_{N\in \mspec{n}}$
(that is, if $f_M$ can 
be written 
as a polynomial in $(a_N)_{N \in S}$ for some finite subset $S$ of $\mspec{n}$),
then $f$ is  a \emph{polynomial map on species}, or simply a 
\emph{polynomial map}.
\end{definition}
\index{polynomial map}

More generally, let $I$ be a finite index set and suppose that $f$ is a 
function from $\prod_{i\in I} \spec{n_i}$ to $\spec{m}$, where the $n_i$ are 
positive integers. Then, as above, $f$ determines functions 
$(f_{M})_{ M \in \mspec{m}}$. If $f_M ( ([R_i])_{i\in I})$
is a polynomial in the coefficients $([R_i])_{i\in I}$,
we will say that $f$ is a polynomial map on species in this
case as well. A map $f:\prod_{i\in I} \spec{n_i} \to \prod_{j\in J} \spec{m_j}$ 
will
be considered a polynomial map on species if each component function is
a polynomial map on species.

\begin{remark}
Any occurrence of $\spec{i}$ in the above definitions may be replaced
by $\specf{i}$, or any other subset of $\spec{i}$; the resulting maps will 
also be considered \emph{polynomial
maps}.
\end{remark}

The chief significance of polynomial maps is that
a polynomial map $f$ on species 
can be extended to a polynomial map on virtual species simply by allowing
the coefficients in the molecular series upon which $f$ operates to take
on negative values. Since two polynomials which agree for 
all positive
values of their arguments must be identically equal, this extension is
the unique polynomial extension of $f$ to virtual species.  Given a polynomial
operation, we will henceforth take this extension for granted and use it without
comment.

We now demonstrate that certain maps on species are polynomial.

\begin{definition}
A binary operation $*$ from $\spec{n} \times \spec{n}$ to 
$\spec{n}$ is  \emph{bilinear} if,
for all $A_i,B_i,C \in \spec{n}$,
\begin{myenumerate}
\item $(A_1+A_2+\dots) * C = A_1*C + A_2*C+\dots$
\item $A*(B_1+B_2+\dots) = A*B_1 + A*B_2+\dots$
\end{myenumerate}
whenever the sums $A_1+A_2+\dots$, $B_1+B_2+\dots$, respectively, are
defined.
\end{definition}

We define \emph{linearity} (of a map from $\spec{n}$ to $\spec{n}$)
entirely analogously.

\begin{definition}
A binary operation $*$ from $\spec{n} \times \spec{n}$ to 
$\spec{n}$ is a \emph{polynomial operation}, or simply \emph{polynomial},
if the map $(A,B)\mapsto A*B$ is a polynomial map on species.
\end{definition}

We note that $A*B$, 
where $A$ and $B$ are virtual species, is well-defined 
for any polynomial binary operation $*$.

If $f$ and $g$ are polynomial maps, then it is immediate that their 
composition, when defined, is also a polynomial map. 

As in \cite{virtspec}, we define a species $A$ to be a 
subspecies \index{subspecies}\index{species!sub-} of the species 
$B$ if
$A[U] \subseteq B[U]$ for all finite sets $U$, and the inclusion is
a natural transformation. If $A$ is a molecular species, this is 
equivalent to saying that $A$ occurs in the molecular series of $B$. 
We note that $*:\spec{n} \times \spec{n} \to \spec{n}$ is
bilinear, then a given molecular species $T$ can be a subspecies
of $M_1 * M_2$ for only finitely many pairs $(M_1,M_2)
\in \mspec{n}\times\mspec{n}$; otherwise,
\begin{equation*}
\biggl( \sum_{M_1\in\mspec{n}} M_1 \biggr) *
\biggl( \sum_{M_2\in\mspec{n}} M_2 \biggr)
\end{equation*}
would not be defined. This observation also gives:

\begin{lemma}
If $*:\spec{n} \times \spec{n} \to 
\spec{n}$ is
bilinear, then it is a polynomial operation.
\label{linkc}
\end{lemma}

\begin{proof}
Let $A=\sum_{M\in \mspec{n}} a_M M$ and $B=\sum_{M\in \mspec{n}} b_M M$
be the molecular series of species $A$ and $B$. Then,

\begin{equation}
A*B = \sum_{M_1,M_2\in \mspec{n}} a_{M_1} b_{M_2} M_1 * M_2
\label{last1}
\end{equation}
by bilinearity. Since a given molecular species $T$ can be a subspecies
of $M_1 * M_2$ for only finitely many pairs $(M_1,M_2)
\in \mspec{n}\times\mspec{n}$,
we can
collect terms on the right-hand side
of \eqref{last1}, yielding a molecular series whose coefficients are
polynomials in the $a_{M_1}$ and $b_{M_2}$.
\end{proof}

\begin{remark}
Similarly, any linear operation is polynomial. Thus, for example,
the diagonal map $\nabla:\spec{2}\to \spec{1}$ 
\index{diagonal map}
(defined
by $(\nabla F)[U] = F[U,U]$), and the map $A\mapsto A'$
from $\spec{1}$ to $\spec{1}$, are polynomial.
\end{remark}

\begin{lemma}
The following binary operations are polynomial: $+$, $\cdot$,
$\times$, $\circ$.
\end{lemma}

\begin{proof}
For $+$ this is immediate, since the coefficient of $M\in \mspec{n}$
in the molecular series of $A+B$ is the sum of the coefficients of
$M$ in the molecular series of $A$ and $B$.
For $\cdot$ and $\times$, which are bilinear, the lemma follows
 immediately from Lemma~\ref{linkc}.
For $\circ$, which is not bilinear, see \cite{virtspec}.
\end{proof}

The Cartesian product in $Y$
\index{Cartesian product!in $Y$}
is clearly bilinear, and thus polynomial.

\begin{lemma}
The map $F(X,Y) \mapsto F(X,1)$ is polynomial.
\label{lem:sub1}
\end{lemma}

\begin{proof}
It is sufficient to observe that given a 2-sorted molecular
species $M(X,Y)$, $M(X,1)$ is a finite sum of 1-sorted molecular species.
\end{proof}

In order to show that the map $\Phi$ introduced in 
Section~\ref{sec:inplethspec} is polynomial, we first show that it has
a certain multiplicative property:

\begin{lemma}
For species $F_1$ and $F_2$,
\begin{equation*}
\Phi (F_1 + F_2) = \Phi (F_1) \times_Y \Phi (F_2)
\end{equation*}
\label{lem:multPhi}
\end{lemma}

\begin{proof}
An element $f \in \Phi (F_1 + F_2)[U,V]$ is a function from
$U$ to the disjoint union $F_1[V] \cup F_2[V]$. Letting $U_i =
f^{-1}(F_i[V])$ for $i=1,2$, we obtain functions $f_i:U_i \to F_i[V]$,
with $U$ the disjoint union of $U_1$ and $U_2$ ($f_i$ is simply $f|_{U_i}$).
This gives a natural bijection from $\Phi (F_1 + F_2)[U,V]$ to
$\sum_{U_1+U_2 = U} \Phi (F_1)[U_1,V] \times \Phi (F_2)[U_2,V]$, proving
the claim. 
\end{proof}

\begin{lemma}
$\Phi : \spec{1} \to \spec{2}$ is a polynomial map on species.
\label{lem:Phi-poly}
\end{lemma}

\begin{proof}
For $A\in\spec{2}$ and $i\in \mathbb{N}$, let
$A^{\times_Y i}$ denote $A\times_Y A \times_Y A \times_Y \dots$
($i$ factors). We take $A^{\times_Y 0}$ to be $1$.

By Lemma~\ref{lem:multPhi}, we have, for any $M\in\mspec{1}$ and
$n\in \mathbb{N}$,
\begin{equation}
\begin{aligned}
\Phi(n M) & = \Phi(M) ^ {\times_Y n} \\
& = (E(Y) + \Phi(M) - E(Y)) \\
& = \sum_{i=0}^n \binom{n}{i} \bigl(\Phi(M) - E(Y)\bigr)^{\times_Y i}
\end{aligned}
\end{equation}
Thus, for any molecular species $M_2\in\mspec{2}$, the coefficient
of $M_2$ in the molecular series of $\Phi(nM)$ is a polynomial
in $n$. Applying Lemma~\ref{lem:multPhi} again, and the
fact that $\times_Y$ is a polynomial operation, we have that
the coefficient of $M_2$ in the molecular series of
$\Phi(a_1 N_1 + a_2 N_2 +\dots)$ is a polynomial in
$a_1,a_2,\dots$, provided that $a_1 N_1 + a_2 N_2 + \dots$
is a finite sum. 

Now consider the coefficient of $M_2\in\mspec{2}$ in the molecular
series of 
\begin{equation*}
\Phi\biggl(\sum_{M\in\mspec{1}} a_M M \biggr)
\end{equation*}
Since $M_2$ is
molecular, it must be homogeneous, say of degree $(n,m)$. By the
definition of $\Phi$, $M_2 [n,m]$ thus consists of the orbit,
under $\sym{n}\times\sym{m}$, of some function 
\begin{equation*}
f:[n]\to 
\biggl(\sum_{M\in\mspec{1}} a_M M\biggr)[m]
\end{equation*}
Thus, we see that the coefficient of $M_2$ depends only upon
those molecular species $M\in\mspec{1}$ which are homogeneous
of degree $m$---in fact, if $S$ is the set of such molecular
species, the coefficient of $M_2$ in $\Phi(\sum_{M\in\mspec{1}} a_M M)$
is equal to the coefficient of $M_2$ in
$\Phi(\sum_{M\in S} a_M M)$. Since this latter is a finite sum,
the lemma is proven.
\end{proof}

\begin{lemma}
Lemma~\ref{lem:multPhi} holds when $F_1$ and $F_2$ are virtual species.
\label{lem:multPhi-virt}
\end{lemma}

\begin{proof}
The maps $(F_1,F_2)\mapsto \Phi(F_1+F_2)$ and $(F_1,F_2)\mapsto 
\Phi (F_1) \times_Y \Phi (F_2)$ are both polynomial maps on species,
and by Lemma~\ref{lem:multPhi}, they agree when $F_1, F_2$ are species.
Therefore, they
agree on virtual species. 
\end{proof}

\begin{lemma}
The map $\inpleth{Y}: \specf{1}\times \spec{2} \to \spec{2}$ 
is a polynomial
map on species.
\end{lemma}

\begin{proof}
Let $F=\sum_{M\in \mspec{1}} a_M M$ and $G=\sum_{N\in \mspec{2}} b_N N$
be the molecular series of $F$ and $G$, and suppose $R$ is a molecular
subspecies of $F\inpleth{Y} G$. We must show that the coefficient of
$R$ in the molecular series of $F\inpleth{Y} G$ is a polynomial in the
$a_M$ and $b_N$.

Since $R$ is molecular, it must be homogeneous of degree $(m_1, m_2)$, for
some $m_1$, $m_2$.  An element of $R$ consists of:
\begin{myenumerate}
\item A partition $U_1 + \dots + U_n = [m_1]$
\item $G$-structures $\alpha_1,\dots,\alpha_n$, with $\alpha_i \in G[U_i,m_2]$
\item An $F$-structure $\alpha \in F[\{\alpha_1,\dots,\alpha_n\}]$
\end{myenumerate}
Thus,
\begin{equation}
\alpha_i \in \sum_{N\in \mspec{2}} b_N N[U_i,m_2]
\label{eq:alpha-where}
\end{equation}
for $1\leq i \leq n$, where
in the sum on the right, $b_N N[U_i,m_2]$ denotes the disjoint union of
$b_N$ copies of the set $N[U_i,m_2]$. (Up to the action of
$\sym{m_1}\times \sym{m_2}$, this is the only element of $R$, since
$R$ is molecular.) 

We will think of the $b_N$ copies of
the set $N[U_i,m_2]$ as each being colored in one of $b_N$ distinct
colors. In~\eqref{eq:alpha-where}, $\alpha_i$ must be an element of
exactly one copy of $N[U_i,m_2]$, for some $N\in \mspec{2}$. We
will say that $\alpha_i$ \emph{occurs} in this copy, and that
$\alpha_i$ occurs in the molecular species $N$.

Let $N_1,\dots,N_s$ be the distinct molecular species in which $\alpha_i$
occurs (for some $i$), and for each $i$, let $k_i$ be the
number of copies of $N_i$ in which some $\alpha_j$ appears. Let $M_1$ be 
the molecular
subspecies of $F$ in which $\alpha$ occurs. Let $a_1 = a_{M_1}$,
and $b_i = b_{N_i}$ for $1\leq i \leq s$. Define
$F_0  = M_1$, $G_0  = \sum_{i=1}^s k_i N_i$.
Then $F_0 \inpleth{Y} G_0$ is naturally isomorphic to a subspecies of
$F\inpleth{Y} G$, and in fact, for each choice of $k_i$ out of $b_i$ copies of
$N_i$, and one copy out of $a_1$ of $M_1$, we obtain a subspecies
of $F\inpleth{Y} G$ which is isomorphic to $F_0 \inpleth{Y} G_0$.
Let $p$ be the coefficient of $R$ in $F_0 \inpleth{Y} G_0$. Then
we obtain the following contribution to the coefficient of $R$ in
the molecular series of $F\inpleth{Y} G$:
\begin{equation}
p a_1 \prod_{i=1}^s \binom{b_i}{k_i} 
\label{eq:Rprod}
\end{equation}

To show that the coefficient of $R$ is in fact a finite sum of such
terms, we proceed as follows. Suppose $R^{*}$ is another molecular subspecies
of $F\inpleth{Y} G$ which is isomorphic to $R$. Then $R^{*}$ must
be homogeneous of degree $(m_1, m_2)$, an element of $R^{*}$ consists
of
\begin{myenumerate}
\item A partition $U_1^{*} + \dots + U_{n^*}^{*} = [m_1]$
\item $G$-structures $\alpha_1^{*},\dots,\alpha_{n^*}^{*}$, with $\alpha_i^{*} \in G[U_i^{*},m_2]$
\item An $F$-structure $\alpha^{*} \in F[\{\alpha_1^{*},\dots,\alpha_{n^*}^{*}\}]$
\end{myenumerate}
and we can define $N_1^{*},\dots,N_{s^*}^{*}$, $k_i^{*}$,
$a_1^{*}$, $b_i^{*}$, $M_1^{*}$, $F_0^{*}$, $G_0^{*}$, and $p^*$ just as we
did their non-starred analogues. If $F_0^{*} \inpleth{Y} G_0^{*} =
F_0 \inpleth{Y} G_0$ (that is, if they are isomorphic), then $R^{*}$ is
counted in the sum~\eqref{eq:Rprod}. If not, we obtain a new contribution
to the coefficient of $R$, equal to:
\begin{equation*}
p^{*} a_1^{*} \prod_{i=1}^{s^*} \binom{b_i^{*}}{k_i^{*}} 
\end{equation*}

There are only finitely many possibilities for the species 
$F_0^{*} \inpleth{Y} G_0^{*}$, which can be seen as follows. Since
$F$ is strictly finite, there are only finitely many possibilities for
$F_0$. And, each
$N_i^{*}$ must be homogeneous of a degree $(d,m_2)$, with
$d\leq m_1$, so there are only finitely many possibilities for the
$N_i^{*}$. Finally, the $k_i^{*}$ are bounded by  coefficients in
the molecular series of $G$.

\end{proof}

\subsection{Polynomial maps on symmetric functions}

A map $f:\Lambda \to \Lambda$ determines a function $f_{\lambda}$
for every partition $\lambda$:

\begin{equation}
f \biggl( \sum_{\lambda \in \partitions} a_{\lambda} \pz{\lambda} \biggr)
= \sum_{\lambda \in \partitions } f_{\lambda} ( (a_{\mu})_{\mu \in \partitions})
  \pz{\lambda}
\label{last2}
\end{equation}

\begin{definition}
If the functions $f_{\lambda}$ in \eqref{last2}  are polynomials in the
$a_{\mu}$, then $f$ is a \emph{polynomial map on symmetric functions},
or simply a \emph{polynomial map}.
\end{definition}

Similarly, a map $f:\spec{m} \to \lam{}{}$ determines functions
$f_{\lambda}$:

\begin{equation}
f (R)
= \sum_{\lambda \in \partitions } f_{\lambda} ( [R] )
  \pz{\lambda} 
\label{last3}
\end{equation}

\begin{definition}
If the functions $f_{\lambda}$ in \eqref{last3}  are polynomials in the
coefficients $[R]$,
then $f$ is a \emph{polynomial map from species to symmetric
functions},
\index{polynomial map}
or simply a \emph{polynomial map}.
\end{definition}

We extend this definition to include polynomial maps from
$\spec{m}\times \spec{n} \to
\lam{}{xy}$, etc., in the obvious way.

\begin{lemma}
The map $F \mapsto Z_F$ is a polynomial map from species to
symmetric functions.
\end{lemma}

\begin{proof}
For simplicity, we consider the case where $F$ is 1-sorted, with molecular series
$F=\sum_{M\in \mspec{1}} a_M M$. Then,
\begin{equation}
Z_F = \sum_{M\in \mspec{n}} a_M Z_M
\end{equation}
The coefficient of $\pz{\lm}$ in $Z_M$, for some partition $\lm$,
can be non-zero for only finitely many $Z_M$ (since there are only
finitely many molecular species homogeneous of a given degree), and
the lemma follows immediately.
\end{proof}

We now have the tools to verify that equalities such as~\eqref{eq:comp-spec}
hold for virtual species. To show that $Z_{F\circ G} = Z_F \circ Z_G$,
we observe that the maps $(F,G)\mapsto Z_{F\circ G}$ and
$(F,G)\mapsto Z_F \circ Z_G$ are both polynomial maps from
species to symmetric functions. By~\eqref{eq:comp-spec}, they agree
on species; therefore, they agree on virtual species.


\chapter{Digraphs}
\label{chap:digraphs}

\section{$G$-digraphs}
\label{sec:G-digraphs}

For a 1-sorted species $G$, we define a $G$-digraph on a
set $U$ to be a
digraph $\mathcal{D}$ with vertex set $U$, together with a
$G$-structure on the set of arcs out of each vertex.
We denote by $\digraphs{G}$ the species of $G$-digraphs in which loops are
not allowed, and by $\digraphsl{G}$ the species of $G$-digraphs in 
which they are. 

\begin{lemma}
$\digraphsl{G} = \nabla \Phi (E \cdot G)$, where $\Phi$ is the map
defined in Section~\ref{sec:inplethspec}.
\label{lem:dig-loop}
\end{lemma}

\begin{proof}
For a finite set $V$, an element of $(E \cdot G)[V]$ consists 
of a subset of $V$, together with a $G$-structure on that subset.
By the definition of $\Phi$, an element of
$\Phi (E \cdot G) [U,V]$  consists of a function from $U$ to 
$(E \cdot G)[V]$---i.e., it associates to each element of $U$
a subset of $V$ and a $G$-structure on that subset.

An element of $\nabla \Phi (E \cdot G)[U]$ thus associates to
 each $u\in U$ a subset of $U$, and a $G$-structure on that subset. This is 
 equivalent to specifying a $G$-digraph on $U$.
\end{proof}

Lemmas~\ref{lem:dig-loop} and \ref{lem:Phi-fix} allow us to calculate
the cycle index of $\digraphsl{G}$ for any species $G$. If
$G=E$, for example, $\digraphsl{G}$ is the species of all digraphs
(with loops allowed), and we find that the isomorphism-types generating
function for $\digraphsl{G}$ is:
\begin{equation*}
\widetilde{\digraphsl{E}} (x) =                                                
2 x + 10 x^2  + 104 x^3  + 3044 x^4  + 291968 x^5  + 96928992 x^6 +\dots 
\end{equation*}
(Counting such digraphs is equivalent to counting relations on 
a set; see Section 5.1 of~\cite{graphenum}.)

\section{Removing loops}
\label{sec:loop-remove}

In order to deal with $D_G$,
we have the following lemma:

\begin{lemma}
If $G$ is a species, then $\digraphs{G+G'} = \digraphsl{G}$.
\label{lem:loop-rid}
\end{lemma}

\begin{proof}
A $G$-structure on the edges out of a vertex with a loop can
also be thought of as a $G'$-structure on these edges with the loop
removed. So given a $\digraphsl{G}$-structure on a set $U$, we can remove 
all the loops
and replace the $G$-structures at those vertices with the same structures
considered as $G'$-structures on one fewer objects. This gives a natural
bijection from  $\digraphsl{G} [U]$ to $\digraphs{G+G'} [U]$. 
\end{proof}

Thus, a way to calculate the cycle index series of
$D_G$ for some species $G$ is to find $G_1$ such that $G_1 +
{G_1}'=G$ and calculate the cycle index series of $\digraphsl{G_1}$. 
We note that one such solution
(provided the sum converges) is $G_1 = G - G' + G'' - G''' + \dots$; the 
$G_1$ so obtained may therefore be a virtual species. We need to extend the
definition of $\digraphsl{G}$ and $\digraphs{G}$ to virtual species and show
that Lemma~\ref{lem:loop-rid} remains true when $G$ is a virtual species.

Lemma~\ref{lem:dig-loop}, combined with 
Lemma~\ref{lem:Phi-poly}, allows us to conclude
that
the map $G \mapsto \digraphsl{G}$
is a polynomial map on species. Therefore $\digraphsl{G}$ is well-defined
for virtual species.

\begin{lemma}
The map $G\mapsto \digraphs{G}$ is a polynomial map on species.
\label{lem:dig-poly}
\end{lemma}

\begin{proof}
Let $M$ be a molecular subspecies of $\digraphs{G}$, where $G=\sum_{N\in
\mspec{1}} a_N N$. Suppose $M$ is homogeneous of degree $m$. An
element of $M[m]$ consists of a digraph $\mathcal{D}$ on $[m]$, together
with a $G$-structure on the edges out of every vertex. Any $G$-structure
is an $N$-structure, for some $N\in \mspec{1}$, colored in one of $a_N$
colors.

Let $N_1,\dots
N_s$, be the elements of $\mspec{1}$ whose structures appear at
the vertices of $\mathcal{D}$. Furthermore,
let $k_i$ be the number of distinct colors of $N_i$-structures which appear,
and let $a_i=a_{N_i}$.

We give an example of such a molecular subspecies in Figure~\ref{graphpic}. 
Here $m=6$, $s=2$, and $G=a_{N_1} N_1 + a_{N_2} N_2 + \dots$ for two 
molecular species
$N_1$ and $N_2$ (homogeneous of degrees 1 and 3, respectively). We see
from the figure that $k_1 =3, k_2 = 1$.

\begin{figure}
\begin{center}
\epsfig{file=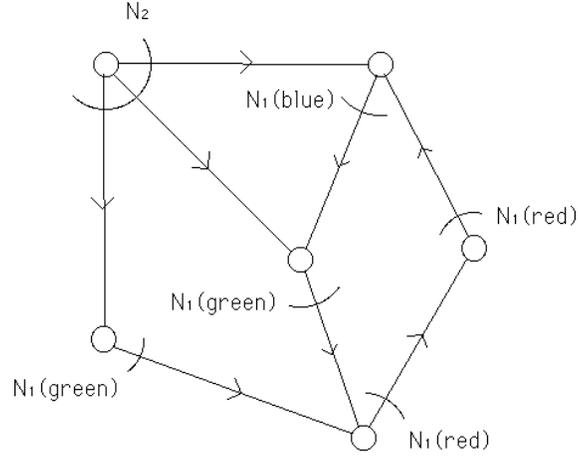,height=2.808in}
\end{center}
\caption{A molecular subspecies}
\label{graphpic}
\end{figure}

Let $T=\sum_{i=1}^{s} k_i N_i$. Then $\digraphs{T}$ is naturally isomorphic to
a subspecies of $\digraphs{G}$, and any choice of $k_i$ colors out of $a_i$ 
yields
a subspecies of $\digraphs{G}$ which is isomorphic to $\digraphs{T}$. 
(Switching red, green,
and blue to yellow, black, and white in Figure~\ref{graphpic}, for example, 
yields
a distinct molecular subspecies isomorphic to the first.) Let $p$ be
the coefficient of $M$ in the molecular series of $\digraphs{T}$. Then we
have the following contribution to the
coefficient of $M$ in the molecular series of $\digraphs{G}$: 
\begin{equation}
p \prod_{i=1}^s \binom{a_i}{k_i} 
\label{last12}
\end{equation}
We will now show that the coefficient of $M$ in $\digraphs{G}$ is in fact
a sum of such terms.

Suppose $M^{*}$ is another molecular subspecies of $\digraphs{G}$ which is
isomorphic to $M$. Then $M^{*}$ must be homogeneous of degree $m$. 
We can define $s^{*}$, $N_1^{*},\dots, N_{s^{*}}^{*}$,
$a_i^{*}$, and $k_i^{*}$ just as we defined their analogues for $M$,
and let $T^{*}=\sum_{i=1}^{s^{*}} k_i^{*} N_i^{*}$. If $T^{*} = T$
then $M^{*}$ is counted in  \eqref{last12}. If not,
let $p^{*}$ be the coefficient of $M^{*}$ in $T^{*}$; we obtain a
new expression of the same form as \eqref{last12}, with all quantities 
replaced by their starred analogues.

To see that there are only finitely many possibilities for the species
$T^{*}$, we note that any molecular species $N_i^{*}$ occurring in the
molecular series of $T^{*}$ must be homogeneous of some degree less
than $m$ (since a vertex in a loopless digraph on $m$ vertices can
have arcs to at most $m-1$ vertices). 
\end{proof}

\begin{lemma}
If $G$ is a virtual species, then $\digraphs{G} = \digraphsl{G+G'}$.
\end{lemma}

\begin{proof}
We recall that the map $A\mapsto A'$ is polynomial (since it is linear).
By our observation that $G \mapsto \digraphsl{G}$ 
is polynomial, we conclude that $G \mapsto \digraphsl{G+G'}$
is polynomial as well. We have just seen that $G \mapsto \digraphs{G}$ is
polynomial, and by Lemma~\ref{lem:loop-rid}, the two maps agree on species.
Therefore, they agree on virtual species.
\end{proof}

It remains to calculate the cycle index series of $\digraphs{G}$ when
$G$ is a virtual species. In order to do this, we will make use of
Lemma~\ref{lem:multPhi-virt}, which shows that for any species $F$, $\Phi (F-F)
=\Phi(F) \times_Y \Phi(-F)$. Since $\Phi (0) = E(Y)$, this gives a
way to calculate the cycle index series of $\Phi (-F)$. In fact, 
applying~\eqref{eq:Phicyc}, we see that
\begin{equation}
Z_{\Phi(-F)} = \sum_{\mu} \exp \biggl( - \sum_{i \geq 1} 
\frac{\fix F[\sigma^i] p_{i}(x)}{i}
\biggr) \frac{p_{\mu}(y)}{z_{\mu}}
\label{eq:Phineg}
\end{equation}
where $\sigma$ is any permutation of cycle type $\mu$.

Combining Lemma~\ref{lem:multPhi-virt} and \eqref{eq:Phineg}
gives a general expression for 
$Z_{\Phi(F_1 - F_2)}$ for any species $F_1$ and $F_2$:
\begin{equation}
\begin{aligned}
Z_{\Phi(F_1 - F_2)} &= \sum_{\mu} \exp \biggl( \sum_{i \geq 1} 
\fix F_1[\sigma^i] \frac{p_{i}(x)}{i}
\biggr) 
\exp \biggl( - \sum_{i \geq 1} 
\fix F_2[\sigma^i] \frac{p_{i}(x)}{i}
\biggr) \frac{p_{\mu}(y)}{z_{\mu}} \\
&= \sum_{\mu} \exp \biggl( \sum_{i \geq 1} 
\bigl(\fix F_1[\sigma^i]-\fix F_2[\sigma^i]\bigr) \frac{p_{i}(x)}{i}
\biggr) \frac{p_{\mu}(y)}{z_{\mu}} \\
&= \sum_{\lambda, \mu} \biggl( \prod_{k \geq 1} 
{\bigl(\fix F_1[\sigma^k]-\fix F_2[\sigma^k]\bigr)}^{\beta_k} \biggr) 
\frac{p_{\lambda}(x)}{z_{\lambda}}
\frac{p_{\mu}(y)}{z_{\mu}}
\end{aligned}
\end{equation}
where $\sigma$ denotes a permutation of cycle type $\mu$, and
$\beta_k$ is the number of $k$-cycles in a permutation of
cycle type $\lambda$.

\section{Applications}

\subsection{All digraphs}

The species of all digraphs is  $\digraphs{E}$, where $E$
is the species of sets.
We saw in Section~\ref{sec:loop-remove} that to calculate
the cycle index of this species, we must find a (possibly virtual)
species $G$ such that $G+G'=E$. Such a species is
$G = 1 + E_2 + E_4 + E_6 + \dots$. Applying~\eqref{eq:Phicyc},
we can calculate the cycle index of $\digraphs{E}$, and
its isomorphism-types generating function:
\begin{equation*}
\begin{aligned}
\widetilde{\digraphs{E}} (x) = \ &                                           
x + 3 x^2  + 16 x^3  + 218 x^4  + 9608 x^5  + 1540944 x^6 
 + 882033440 x^7 \\
& + 1793359192848 x^8 + 13027956824399552 x^9 +\dots
\end{aligned}
\end{equation*}

\subsection{Digraphs in which every vertex has outdegree $k$}

The species of such digraphs is $\digraphs{E_k}$. So we must solve
$G+G'=E_k$ for $G$. The virtual species $E_k - E_{k-1} + E_{k-2} +\dots$
is such a solution. We summarize our results in Table~\ref{outtable}.

\begin{table}
\begin{tabular}{c|r|r|r|r|r|}
 \multicolumn{1}{c}{} & \multicolumn{5}{c}{Outdegree} \\ 
\mcc{$n$} & \mcc{1} & \mcc{2} & \mcc{3} & \mcc{4} & \mcc{5} \\ \hline
2 & 1 & 0 & 0 & 0 & 0 \\
3 & 2 & 1 & 0 & 0 & 0 \\
4 & 6 & 6 & 1 & 0 & 0 \\
5 & 13 & 79 & 13 & 1 & 0 \\
6 & 40 & 1499 & 1499 & 40 & 1 \\
7 & 100 & 35317 & 257290 & 35317 & 100 \\
8 & 291 & 967255 & 56150820 & 56150820 & 967255 \\
9 & 797 & 29949217 & 14971125930 & 111359017198 & 14971125930
\end{tabular}
\caption{Digraphs of outdegree $k$ on $n$ vertices}
\label{outtable}
\end{table}

\subsection{Digraphs with outdegrees from a prescribed set}

Given a specified set $S$ of positive integers, we can enumerate digraphs
in which all outdegrees are members of $S$. Consider the case
$S= \{ 1,3,4 \}$, for example. The species of digraphs in which all
outdegrees are members of $S$ is $\digraphs{G}$, where
$G = E_4 + E_3 + E_1$. A solution to $G_1 + G_1' = E_4 + E_3 + E_1$
is $G_1=E_4 + E_1 -1$, and we obtain the following isomorphism-types
generating function:
\begin{equation*}
\widetilde{\digraphs{G}} =  
     x^2  + 2 x^3  + 19 x^4 + 616 x^5  + 93815 x^6  + 39097411 x^7
  + 30749550146 x^8 \dots
\end{equation*}

\chapter{Graphs}
\label{chap:graphs}

\section{$G$-graphs}

In analogy with Section~\ref{sec:G-digraphs}, we define a
$G$-graph on a
set $U$ to be a
graph  with vertex set $U$, together with a
$G$-structure on the set of vertices incident at each vertex.
We denote by $\graphsl{G}$ the species of $G$-graphs in which loops 
and multiple edges are
allowed, and by $\graphs{G}$ the species of $G$-graphs in which
loops are not allowed (but multiple edges are). 

We note that since multiple edges are allowed, 
$\graphs{G}$ and $\graphsl{G}$ are defined \myiff\ $G$
is strictly finite.

\begin{figure}
\begin{center}
\epsfig{file=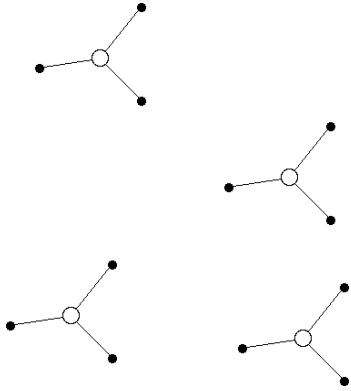,height=2in}
\end{center}
\caption{An $E(X \cdot E_3 (Y))$-structure}
\label{fig:exeky}
\end{figure}

\begin{figure}
\begin{center}
\epsfig{file=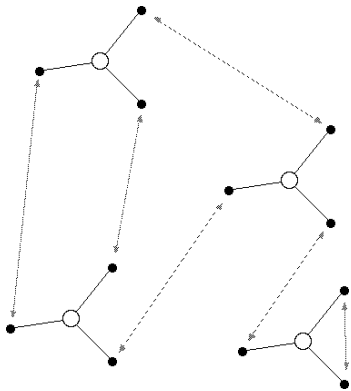,height=2in}
\end{center}
\caption{A 3-regular graph}
\label{fig:3reg}
\end{figure}

To construct the species of $G$-graphs combinatorially, consider
first the species 
$X \cdot G (Y)$ and $E(X \cdot G(Y))$. A structure of this latter species
is pictured in Figure~\ref{fig:exeky} with $G=E_3$. Creating a $G$-graph
from such
a structure amounts to pairing up the $Y$ points to form edges, as shown
in Figure~\ref{fig:3reg}.
(We note that both loops and multiple edges are
permitted by this construction.)

In species-theoretic terms this means that  a $G$-graph
is specified by both an $E(X \cdot G (Y))$-structure
and an $E(E_2(Y))$-structure on a
given set of points. The species of such pairs of structures
is $E(X \cdot G (Y)) \times_Y E(E_2(Y))$. Since each edge consists
of \emph{two} $Y$-points, the number of $Y$ points is twice the number
of edges in the corresponding graph. In any case, setting $Y=1$ gives the 
species
$\graphsl{G}$:
\begin{equation}
\begin{aligned}
\graphsl{G} (X) & = \bigl( E(X \cdot G (Y)) \times_Y E(E_2(Y)) \bigr) \big|_{Y=1} \\
& = \bigl\langle E(X \cdot G (Y)) , E(E_2(Y)) \bigr\rangle_Y
\end{aligned}
\label{eq:graphsldef}
\end{equation}
Recalling Lemma~\ref{lem:sub1}, and observing that the operations involved
are polynomial,
we note that the map $G\mapsto \graphsl{G}$ is polynomial.

\section{Removing loops}
\label{sec:graphs-loop-remove}

Removing loops from $G$-graphs is  analogous to removing
them from $G$-digraphs. For $n>0$, let $G^{(n)}$ denote the
result of applying the derivative operator to $G$ $n$ times.
We have:

\begin{lemma}
If $G$ is a species, then $\graphs{G+G''+G^{(4)}+G^{(6)}+\dots} = \graphsl{G}$.
\label{lem:graphs-loop-rid}
\end{lemma}

\begin{proof}
A $G$-structure on the edges incident upon a vertex with a loop can
also be thought of as a $G''$-structure on these edges with the loop
removed. If two loops are present, we can regard a $G$-structure
on the incident edges as a $G^{(4)}$-structure on these edges with
both loops removed, and so forth.
Thus, given a $\graphsl{G}$-structure on a set $U$, we can remove 
all the loops
and replace the $G$-structures at those vertices with the same structures
considered as $G''$-structures on two fewer objects, or
$G^{(4)}$-structures on four fewer objects, etc.
This gives a natural
bijection from  $\graphsl{G} [U]$ to $\graphs{G+G''+G^{(4)}+G^{(6)}+\dots} [U]$. 
\end{proof}

\begin{lemma}
The map $G\mapsto \graphs{G}$ is polynomial.
\label{lem:graph-poly}
\end{lemma}

\begin{proof}
The proof is  analogous to that of Lemma~\ref{lem:dig-poly}.
Let $M$ be a molecular subspecies of $\graphs{G}$, where $G=\sum_{N\in
\mspec{1}} a_N N$. Suppose $M$ is homogeneous of degree $m$. An
element of $M[m]$ consists of a graph $\mathcal{R}$ on $[m]$, together
with a $G$-structure on the edges incident upon every vertex. Any $G$-structure
is an $N$-structure, for some $N\in \mspec{1}$, colored in one of $a_N$
colors.

Let $N_1,\dots
N_s$, be the elements of $\mspec{1}$ whose structures appear at
the vertices of $\mathcal{R}$. Furthermore,
let $k_i$ be the number of distinct colors of $N_i$-structures which appear,
and let $a_i=a_{N_i}$.

Let $T=\sum_{i=1}^{s} k_i N_i$. Then $\graphs{T}$ is naturally isomorphic to
a subspecies of $\graphs{G}$, and any choice of $k_i$ colors out of $a_i$ 
yields
a subspecies of $\graphs{G}$ which is isomorphic to $\graphs{T}$. 
 Let $p$ be
the coefficient of $M$ in the molecular series of $\graphs{T}$. Then we
have the following contribution to the
coefficient of $M$ in the molecular series of $\graphs{G}$: 
\begin{equation}
p \prod_{i=1}^s \binom{a_i}{k_i} 
\label{eq:graphproof-last}
\end{equation}
We will now show that the coefficient of $M$ in $\graphs{G}$ is in fact
a sum of such terms.

Suppose $M^{*}$ is another molecular subspecies of $\graphs{G}$ which is
isomorphic to $M$. Then $M^{*}$ must be homogeneous of degree $m$. 
We can define $s^{*}$, $N_1^{*},\dots, N_{s^{*}}^{*}$,
$a_i^{*}$, and $k_i^{*}$ just as we defined their analogues for $M$,
and let $T^{*}=\sum_{i=1}^{s^{*}} k_i^{*} N_i^{*}$. If $T^{*} = T$
then $M^{*}$ is counted in  \eqref{eq:graphproof-last}. If not,
let $p^{*}$ be the coefficient of $M^{*}$ in $T^{*}$; we obtain a
new expression of the same form as \eqref{eq:graphproof-last}, 
with all quantities 
replaced by their starred analogues.

To see that there are only finitely many possibilities for the species
$T^{*}$, we use the fact that $G$ is strictly finite (recall that
$\graphs{G}$ is only defined if $G$ is strictly finite), and therefore
has only finitely many molecular subspecies.
\end{proof}

The map $G\mapsto \graphsl{G}$ is polynomial.  Thus, Lemma~\ref{lem:graph-poly}
allows us to conclude that
Lemma~\ref{lem:graphs-loop-rid} holds for virtual species. To
calculate the cycle index of $\graphs{G}$, we can solve
the equation $G_1 + G_1'' + {G_1}^{(4)} + \dots = G$ for $G_1$, 
and calculate the cycle
index of $\graphsl{G_1}$. We note that  $G - G''$ is a solution.

\section{Application: regular graphs}

Consider the problem of enumerating $k$-regular graphs, that is, graphs
in which every vertex has degree $k$.
In our framework, these graphs are
given by $\graphs{E_k}$ or $\graphsl{E_k}$, depending on whether or not
loops are allowed. 

Consider $\graphs{E_3}$, for example. Applying
the methods of Section~\ref{sec:graphs-loop-remove}, we solve
$G_1 + G_1'' + {G_1}^{(4)} + \dots = E_3$
for $G_1$, obtaining $G_1 = E_3 - E_1$.
We then use the combinatorial construction~\eqref{eq:graphsldef} to
calculate the cycle index of $\graphsl{E_3 - E_1}$, which, by
Lemma~\ref{lem:graphs-loop-rid}, is also the cycle index of
$\graphs{E_3}$. We obtain the following isomorphism-types generating 
function:
\begin{equation}
\widetilde{\graphs{E_3}}(x) =
x^2 + 3x^4 + 9x^6 + 32x^8 + 135x^{10} + \dots
\end{equation}

We note that the problem of enumerating 3-regular graphs arose in
chemistry in 1966 (though the interest there was in connected graphs).
Read~\cite{read1} had solved 
the problem of enumerating
$3$- and $4$-regular graphs in which loops, multiple edges, or both, are excluded,
but only in the case in which either the edges or the vertices were unlabeled.
The general problem was considered unsolved.

In fact, the Dutch mathematician Jan de Vries had enumerated
$k$-regular graphs in the unlabeled case in 1891 
(though not in graph-theoretic language), but his work did not
become widely known for nearly a century. See~\cite{enum100} for
the full story.

Asymptotic results are also known (see \cite{asymp}).

\section{Bicolored graphs}

We recall that a \emph{bicolored graph} is one whose vertices have been 
partitioned into two (non-empty) sets, such that vertices in the same set 
are never
adjacent. We define a \emph{bicolored $G$-graph} to be a $G$-graph
in which the vertices have been so partitioned.
To count  bicolored $G$-graphs, we use the inner plethysm
in $Y$.

\begin{figure}
\begin{center}
\epsfig{file=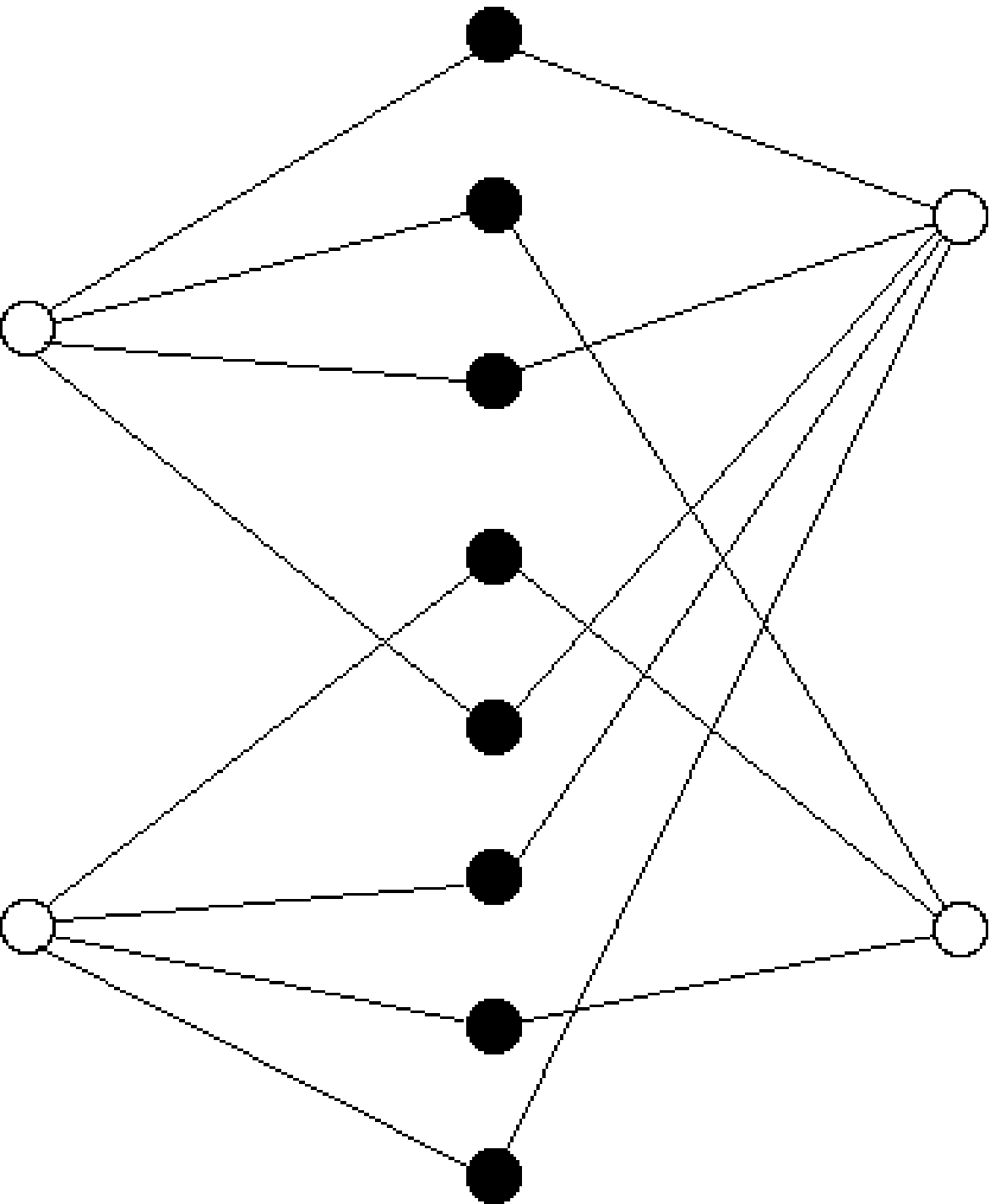,height=1.713in}
\end{center}
\caption{An $E_2 \inpleth{Y} E^{*}(X \cdot E(Y))$-structure}
\label{fig:bi-graph}
\end{figure}

Let $E^{*} = E_1+E_2+E_3+\dots$ denote the species of non-empty sets,
and consider the species
$E^{*}(X \cdot G(Y))$.
A structure
of $E_2 \inpleth{Y} E^{*}(X \cdot G(Y))$ consists of two
such structures which share the same $Y$ points 
(see Figure~\ref{fig:bi-graph} for an example when
$G=E$). Thus, $E_2 \inpleth{Y} E^{*}(X \cdot G(Y))
=\bigraphs{G}(X,Y)$, where $\bigraphs{G}(X,Y)$ is the
species of bicolored $G$-graphs with vertices
of sort $X$ and edges of sort $Y$.

The techniques of Section~\ref{sec:inplethY-construction} allow us to
calculate the cycle index series of $E_2 \inpleth{Y} E^{*}(X \cdot G(Y))$,
and thus its isomorphism-types generating function.  For example,
with $G=E$, we find:
\begin{equation*}
\begin{aligned}
\widetilde{\bigraphs{E}}(x,y) = \ &
x^2 (1+y+y^2+y^3+y^4+y^5+\dots) \\
&+x^3 (1+y+2y^2+2y^3+3y^4+3y^5+\dots) \\
&+ x^4 (2 + 2y + 5y^2 + 7y^3 + 12y^4 + 15y^5+\dots) \\
&+ x^5 (2 + 2y + 6y^2 + 10y^3 + 21y^4 + 32y^5+\dots) \\
&+ \dots
\end{aligned}
\end{equation*}

We note that the problem of enumerating bicolored graphs, in which
multiple edges are not allowed, was solved by Parthasarathy 
in~\cite{givenpart}.

\renewcommand{\baselinestretch}{1}
\tiny\normalsize

\bibliographystyle{amsplain}
\bibliography{GS}

\printindex

\end{document}